\documentclass[a4paper]{article}

\usepackage[utf8]{inputenc}
\usepackage{fontspec}
\usepackage[english]{babel}
\usepackage{csquotes}

\usepackage{amsmath, amssymb, amsthm, mathtools}
\usepackage{dsfont}        

\usepackage{graphicx}
\usepackage{geometry}
\usepackage{caption}
\usepackage{subcaption}
\usepackage{float}
\usepackage{multirow}
\usepackage{array}
\usepackage{booktabs}
\usepackage{hyperref}
\usepackage[style=numeric,sorting=none]{biblatex}
\usepackage{icomma}        
\usepackage{color}         
\usepackage{todonotes}     
\usepackage{changes}

\newcommand{\med}{\mathrm{med}}
\newcommand{\indicator}{{\mathds 1}}
\newcommand{\id}{\mathrm{d}} 
\newcommand{\e}{\mathrm{e}}
\newcommand{\supp}{\mathrm{supp}\,} 
\newcommand{\BE}{{\mathbb{E}}}

\newcommand{\BR}{{\mathbb{R}}}
\newcommand{\BZ}{{\mathbb{Z}}}
\newcommand{\CB}{{\mathcal{B}}}
\newcommand{\CC}{{\mathcal{C}}}
\newcommand{\CF}{{\mathcal{F}}}
\newcommand{\CM}{{\mathcal{M}}}
\newcommand{\CO}{{\mathcal{O}}}
\newcommand{\SFC}{{\mathsf{C}}}

\theoremstyle{plain}
\newtheorem{theorem}{Theorem}[section]

\newtheorem{lemma}{Lemma}[section]
\newtheorem{corollary}{Corollary}[section]

\theoremstyle{definition}
\newtheorem{definition}{Definition}[section]
\newtheorem{problem}{Problem}
\newtheorem{question}{Question}
\newtheorem{example}{Example}[section]
\newtheorem{remark}{Remark}[section]

\providecommand{\keywords}[1]
{
  \small	
  \textbf{\textit{Keywords---}} #1
}

\begin{document}
\pagenumbering{arabic}
\title{Quantization of Probability Distributions via Divide-and-Conquer: Convergence and Error Propagation under Distributional Arithmetic Operations}
\author{
Bilgesu Arif Bilgin
\thanks{
Signaloid
},  
Olof Hallqvist Elias
\thanks{
Signaloid
}, 
Michael Selby
\thanks{
Signaloid
}, 
Phillip Stanley-Marbell
\thanks{
Signaloid
}
}
\date{\today}

\maketitle
\begin{abstract}
This article studies a general divide-and-conquer algorithm for approximating continuous one-dimensional probability distributions with finite mean.
The article presents a numerical study that compares pre-existing approximation schemes with a special focus on the stability of the discrete approximations when they undergo arithmetic operations.
The main results are a simple upper bound of the approximation error in terms of the Wasserstein-1 distance that is valid for all continuous distributions with finite mean. In many use-cases, the studied method achieves optimal rate of convergence, and numerical experiments show that the algorithm is more stable than pre-existing approximation schemes in the context of arithmetic operations.
\end{abstract}
\section{Introduction}\label{s.intro}
Computers carry out arithmetic operations with point-valued numbers.
The data that dominate contemporary computing systems are however from measurement processes such as sensors.
All measurements are inherently uncertain and this uncertainty is often characterized statistically and forms aleatoric uncertainty.
In addition, many other contemporary applications of probability distributions, such as machine learning, comprise models which also have inherent epistemic uncertainty (e.g., on the weights of a neural network \cite{pmlr-v37-blundell15,NIPS2017_2650d608,KIUREGHIAN2009105}). Hardware and software can exploit this uncertainty in measurements for improved performance as well as for trading performance for power dissipation or quality of results \cite{crayon,stanley2006programming,Stanley_Marbell_2016}.  At the same time, applications such as the numerical solution of stochastic differential equations (SDEs) can be easier to implement and can run faster, when hardware can natively perform arithmetics on discrete representations of probability distributions. All of these applications stand to benefit from more effective methods for representing arbitrary real-world probability distributions and for propagating those through distributional arithmetic operations. A distributional arithmetic operation is defined as the binary operation on probability distributions induced by applying the corresponding arithmetic operation to random variables governed by those distributions.

Recent research has investigated methods to efficiently represent probability distributions using a variety of approaches. Among these, several approaches are based on the method of moments \cite{MR4564849,OSOGAMI2006524,MR4134783} while others formulate the representation problem as an optimization problem for a given metric \cite{MR4647975}. In contrast, the algorithm introduced in this article does not rely on solving an optimization problem or on moment-matching, instead using a recursive domain-splitting procedure that requires only the ability to evaluate a summary statistic (such as the mean or median) of conditional distributions. In addition, ad hoc software has been developed to automate handling uncertainty tracking in computations \cite{Bornholt_2014,giordano2016uncertaintypropagationfunctionallycorrelated,844075,Lafarge01092015} and methods to bound the distributions of errors as they propagate through arithmetic operations \cite{Uncertainty-Propagation-Using-Probabilistic,JSSv076i01,MR2970251,sound-probabilistic-numerical}.
While there has been work on efficient machine representations for interval-valued numbers in computing architectures \cite{standard-for-interval-arithmetic,intro-to-standard-interval-arithmetic} and for approximate real numbers \cite{standard-for-floating-point}, the challenge of efficient representations for uncertain values in computers remains an open problem.

In many applications, the probability distribution one seeks to approximate is often not available in closed form and can only be evaluated through a long sequence of arithmetic/functional operations. Thus, one needs to resort to some form of approximation scheme. The Monte Carlo method currently provides the most popular method that allows for the approximation of distributions whose distribution or density functions (CDFs or PDFs) are not available in closed form. However, the rate of convergence is relatively slow: Under the necessary and sufficient assumption, \cite[Theorem 1.1, Theorem 2.1]{MR1698999}, on the cumulative distribution function $F$ of the true distribution
\[
\int_\BR \sqrt{F(t)(1-F(t))} \id t < \infty,
\]
the expected Wasserstein-1 distance between the empirical and true distribution decays like $1/\sqrt{N}$, where $N$ is the number of independent samples obtained from the distribution $F$. The stochastic nature of Monte Carlo introduces inherent variability, making it difficult in practice to determine how large $N$ must be to achieve a desired level of approximation fidelity. Monte Carlo methods typically involve generating a large number of independent samples from some collection of distributions (the empirical input measures), which are then aggregated through a sequence of operations or transformations to produce an empirical measure of the output distribution (the empirical output measure). At the same time however, it is unclear how the closeness between the input empirical measures and their true distributions translates to the closeness between the output empirical measure and its true distribution. Moreover, due to the inherent randomness of Monte Carlo sampling, the approximation fidelity of the input empirical measures will vary. Given the assumed independence, the probability that some of these input empirical measures fail to accurately approximate their true distributions is non-negligible. For example, suppose we conduct a Monte Carlo simulation with $n$ independent random variables as input and that the probability of the event that the Wasserstein-1 distance between the $i$-th empirical input measure and its true distribution is smaller than $\delta>0$ equals $1-\epsilon,\ \epsilon >0$ for some small $\epsilon,\ \delta$ independent of $n$. Then, the probability that at least one of the empirical input distributions fails to approximate their true distribution with an error of $\delta$ equals 
\[
1-(1-\epsilon)^n\geq 1-\e^{\epsilon n}.
\]
Hence, for a fixed $\epsilon$, this probability tends to $1$ as the number of independent input random variables tends to infinity.

Finding the optimal representation with respect to a certain metric is in general non-trivial and one typically relies on optimization algorithms that seek to minimize the metric of choice. Along with this follow the usual potential problems of optimization such as, potential lack of convexity, slow convergence, ill conditioning and numerical instability, \cite{MR2244940,MR2061575}, making this approach infeasible in systems where approximating distributions needs to be automated and consistently accurate. It is worth noting that when the metric is chosen from the family of Wasserstein metrics, the problem of finding the optimal representation falls into the classical field of \emph{quantization of probability measures} (see the monographs \cite{MR1764176,MR4704055} and references therein).

A few notable exceptions where there are explicit optimal representations are for the Laplace and Exponential distributions with respect to the Wasserstein-1 metric, see \cite[Examples 5.6,5.7 on p.69-70]{MR1764176} and the optimal representation with respect to the Cramer-von Misés distance, \cite{kennan2006note,MR4647975}.
There are, on the other hand, \emph{asymptotically-optimal} representations, \cite[Section 7 p.93]{MR1764176} for an outline of the approach, but the numerical experiments conducted in this paper seem to indicate that in many cases the algorithm introduced in this article performs just as well as the asymptotically-optimal approach with regards to arithmetic operations.
Moreover, the asymptotically optimal representation requires considerably more assumptions on the distribution,  \cite{MR730904,linder1991asymptotically,MR1275992}.

The algorithm introduced in this article is performant with \emph{minimal} assumptions on the distribution, while still achieving close to optimal accuracy.

\subsection{Contributions}
This article makes the following contributions to the theory and numerical evaluation of discrete representations of continuous probability distributions, with particular attention to how approximation error propagates under distributional arithmetic operations.

The article introduces a divide-and-conquer domain-splitting algorithm (that depends on a user-specified rule for splitting the domain) that computes discrete approximations (or representations) of continuous probability distributions possessing finite mean. A special case of this algorithm was considered in \cite{9756254} under the name \emph{Telescoping Torques Representation} or TTR for short.

One key benefit with the algorithm compared to pre-existing ones, as Section~\ref{s.numerical} and~\ref{s.discussion} show, is that certain statistical quantities (depending on the chosen domain-splitting rule) are stable under arithmetic operations. Compared with previous methods of approximating distributions, \cite{MR730904,linder1991asymptotically,MR1275992,MR1764176}, the algorithm requires minimal assumptions on the distribution.

The algorithm exhibits, in many cases, the optimal rate of decay of the fidelity of its approximation when measured with the Wasserstein-1 distance.
Here, \emph{optimal rate of decay} refers to the Wasserstein-1 distance between the target measure and its approximation decreasing at a rate of \((1 + o(1))/N\), where \(N\) denotes the representation size. This result is commonly referred to as Zador's theorem which was originally proved for the $L^2$-norm, \cite{MR2614227}, but a general $L^p$ version is available \cite[Theorem~6.2,p.~78]{MR1764176}.

The article presents a numerical study comparing different representations and examining how approximation error propagates when performing arithmetic operations between distributions. We benchmark the newly implemented representations against Monte Carlo simulations by comparing the expected rate of convergence obtained in \cite{MR1698999} and the deterministic error obtained from our representations.
\subsection{Theoretical main result}
We shall state our main result for a special case of the algorithm. For the general result, see Theorem~\ref{thm.unbounded-support-error-bound}.

The algorithm is a domain-splitting, divide-and-conquer, recursive algorithm that takes a \emph{continuous} probability distribution $\mu$ and a number $n\geq0$ as input and returns a discrete probability distribution of $2^n$ distinct Dirac measures.
Let $T(n,\mu)$ be the probability distribution generated by said algorithm, then we may compute $T(n,\mu)$ via the following steps:
\begin{enumerate}
    \item If $n = 0$, return a Dirac measure at the mean $\bar{\mu}$: $T(\mu,0) = \delta_{\bar{\mu}}$.
    \item If $n\geq 1$, decrement $n$, split the support of $\mu$, which we denote by $\Omega$, into two pieces $\Omega_- = \{ x : x < \bar{\mu}\}$ and $\Omega_+ = \Omega_-^c$.
    Apply the algorithm to the conditional distributions of $\mu$ on $\Omega_\pm$ respectively, i.e., $\mu_\pm := \mu( \cdot \cap \Omega_\pm)/\mu(\Omega_\pm)$ and weight the result by the original mass $\mu(\Omega_\pm)$:
    \begin{equation*}
    T(\mu,n) = \mu(\Omega_-) T(\mu_-,n-1)+\mu(\Omega_+) T(\mu_+,n-1).
    \end{equation*}
\end{enumerate}
Our main result in this case is:
\begin{theorem}
    Assume that $\mu$ is a continuous probability measure with finite mean and assume that $\supp(\mu) = \BR_+$ and let $\mu^{(n)} = T(\mu,n)$. Let $X$ be a random variable with distribution $\mu$, $\omega_{-1}= 0,\ \omega_0 = \int_{\BR_+} t \id \mu(t)$ and define $\omega_{j+1} = \BE[X | X \geq \omega_{j}], j\geq 0$ and $\Omega_j :=[\omega_{j-1},\omega_{j}]$, then for any $n\geq 0$ we have
    \begin{equation*}
        \begin{split}        
            W_1(\mu,\mu^{(n)}) 
            \leq
            \sum_{j=0}^{n-1} 
                \frac{
                        (\omega_{j}-\omega_{j-1})\mu(\Omega_{j})
                    }
                    {2^{n-j}} 
            + \BE[ |X - \omega_{n}| ; X \geq \omega_{n-1}].
        \end{split}
    \end{equation*}
\end{theorem}

We now briefly explain how the result above is obtained. As a first step one proves that for any probability measure $\nu$ with support $\supp(\nu) \subseteq [a,b]$ for some $-\infty < a < b<\infty$ we have
\begin{equation*}
W_1(\nu,T(\nu,n))\leq \frac{b-a}{2^{n+1}},
\end{equation*}
see Theorem~\ref{thm.compact-mean-bound}.
One then uses this result to bootstrap the argument by exploiting the recursive nature of the algorithm, see Corollary~\ref{c.recursive-identity}.

As a corollary, we show that for distributions with polynomially-decaying tails, that is, distributions satisfying $\lim_{x\to \infty} x^\alpha\mu[x,\infty) = c \in (0,\infty)$ where $\alpha>1$, we essentially get the optimal decay rate whenever $\alpha>2$.
More precisely:
\begin{theorem}\label{thm.poly-tails}
    Assume that $\lim_{x\to \infty} x^\alpha\mu[x,\infty) = c \in (0,\infty),\ \alpha>1$ as $x \to \infty$ and let $\mu^{(n)} = T(\mu,n)$. Then we have
    \begin{equation*}
        \lim_{n \to \infty} \frac{\log(W_1(\mu,\mu^{(n)}))}{n} = 
        \log 
        \left(
        \left(
            1-\frac{1}{\alpha}
        \right)^{(\alpha-1)}
        \vee
        \frac{1}{2}
        \right).
    \end{equation*}
\end{theorem}

\subsection{Outline of the article}
Section~\ref{s.not-def} introduces the necessary concepts such as the Wasserstein distance and Section~\ref{s.algorithm} and the algorithm in question. Section~\ref{s.analysis} contains the mathematical analysis of the algorithm and the proofs of the main results. Section~\ref{s.numerical} presents a numerical study focusing on the stability of the approximations with respect to arithmetic operations and Section~\ref{s.discussion} summarizes the findings of this article.

\section{Notation, definitions and prerequisites}\label{s.not-def}
The article uses the standard "big O" notation.
Given two functions $f,\ g :X \to \BR$ we write $f(x) = \CO(g(x))$ as $x \to x_0$ if 
$\limsup_{x \to x_0} f(x)/g(x)< \infty$ while if $\lim_{x \to x_0} f(x)/g(x) = 0$ we write 
$f(x) = o(g(x))$. 
Moreover, we say, $f(x) = \Theta(g(x))$ as $x \to x_0$ if $0<\liminf_{x \to x_0} f(x)/g(x) \leq \limsup_{x\to x_0} f(x)/g(x)<\infty$. Finally, we define $f \sim g, x \to x_0$ if $\lim_{x \to x_0} f(x)/g(x)  = 1.$

Let $(\Omega,\CF)$ denote a measurable space. Given a set $A$ we define the indicator function of $A \subset \Omega$ by $\indicator_A$ or $\indicator\{A\}$ and we let $\#(A)$ denote its cardinality. If $A \subset \CB(\BR)$, where $\CB(\BR)$ denotes the Borel-sigma algebra, we let $|A|$ denote the Lebesgue measure of the set.

For a probability measure $\mu$ on $(\BR,\CB(\BR))$ we let 
\begin{equation*}
    F_\mu(x) := \mu(-\infty,x], \ x \in \BR,
\end{equation*}
denote the cumulative distribution function (CDF). Given two probability measures $\mu,\nu$ on $(\BR,\CB(\BR))$ we let 
\[
\CF \ni A \mapsto (\mu*\nu)(A),
\] denote the convolution of the two measures. Moreover we define the median, $\med(\mu)$, as the smallest value $m$ such that 
\begin{equation}
\mu(-\infty, m-] \leq \frac{1}{2} \leq \mu(-\infty,m],
\end{equation}
where $\mu(-\infty, m-] := \lim_{x \uparrow m} \mu(-\infty, x],$ and the mean $\bar{\mu} := \int_{\BR} x \mu(\id x). $    

Throughout the paper we shall use the following notation and so we introduce it here for clarity. By the notation $\alpha \in \{\pm\}^n$ we mean a string of length $n$ consisting of $+$'s and $-$'s with the convention $\{\pm\}^0=\emptyset$. Moreover, given two strings $\alpha \in \{\pm\}^n, \beta \in \{\pm\}^m$ we define the concatenation 
\(
\alpha \beta \in \{\pm\}^{n+m}
\)
by
\[
\alpha \beta:=\alpha_1\dots\alpha_n\beta_1\dots\beta_m,
\]
with the convention 
\(\emptyset\beta = \beta.\)
\subsection{Wasserstein spaces}
Given a measurable space $(\Omega, \CF)$ we let $\CM(\Omega, \CF)$ denote the space of probability measures and when it is clear from context we will drop the emphasis on the sigma-algebra and just write $\CM(\Omega)$ or $\CM$.
Moreover, we let $\CM_r \subset \CM, r>0$ denote the subspace consisting of all probability measures with finite $r$'th moment.  

For any two probability measures, $\mu,\ \nu \in \CM_r$, $r\geq1$ on a metric space $(X,d)$ we define the Wasserstein metric \cite[Definition 6.1, p.105]{MR2459454} as
\begin{equation*}
    W_p(\mu,\nu) =\left( \inf \int_{X \times X} d(x,y)^p \id \gamma(x,y) \right)^{1/p},
\end{equation*}  
where the infimum is over all probability measures on $X \times X$ satisfying $\gamma(A \times X ) =\mu(A), \gamma(X \times A ) = \nu(A),$ or in the probabilistic language over all couplings between $\mu$ and $\nu$.
Moreover, if $X=\BR$ then it is known, \cite[Theorem 3.1.2,p 109]{MR1619170} and \cite[Proposition 2.17, p.66]{MR3409718}, that
\begin{equation*}
    W_1(\mu,\nu) = \int_\BR \left | F_\mu(x)-F_\nu(x)\right | \id x.
\end{equation*}
While we shall not need them, we note that for $p>1$ there is an analogous formula for $W_p$.
It is a well-known fact that Wasserstein distances metrize weak convergence, \cite[Theorem 6.9]{MR2459454}.

We now discuss how one can determine the lower bound of the Wasserstein-1 error and then how to explicitly compute optimal and asymptotically optimal representations in one dimension.

\subsection{Zador's Theorem}
First we recall Zador's theorem, \cite[Theorem 6.2, p.78]{MR1764176}, that determines the theoretical lower bound of the Wasserstein-1 error for a given discrete representation size. We will use this bound as both a theoretical and practical benchmark against the algorithms discussed in this article.
\begin{theorem}[Theorem~6.2, p.~78 in \cite{MR1764176} ]\label{thm.zador}
    Suppose $\BE[|X|^{1+\delta}]<\infty$ for some $\delta>0$ and assume that the distribution of $X$, denoted by $\mu$, is absolutely continuous with respect to Lebesgue measure with density $f$. 
    Then 
    \begin{equation*}
        \lim_{n\to \infty} n \inf_{\substack{\nu \in \CM \\ \#(\supp(\nu)) =n }} W_{1}(\mu,\nu) = \frac{1}{4} 
        \left(
            \int \sqrt{f} \id x 
        \right)^2,
    \end{equation*}
    where the infimum is over all probability measures with $n$ distinct Dirac measures.
\end{theorem}
\begin{remark}
    This remarkable result provides a lower bound for the Wasserstein-1 error:
    \[
    W_{1}(\mu,\mu^{(n)}) \geq \frac{c(\mu)}{2^n}.
    \]
\end{remark}    

\subsection{Explicit optimal quantizer for Wasserstein-1}
If one wishes to explicitly compute an optimal quantization with respect to $W_1$, then this is achieved by solving a system of equations \cite[Section 5.2]{MR1764176}. We remark that uniqueness is guaranteed only when the distribution is strongly unimodal \cite[Theorem 5.1]{MR1764176}.

If 
\begin{equation*}
    \mu_{\text{opt},n} = \sum_{i=1}^n p_i \delta_{x_i }
\end{equation*} denotes the optimal representation with representation size $n$, then $p_i,\ x_i, i = 1,2,\dots,n$ are given by the following. Let 
\(m_i = (x_i+x_{i+1})/2,\ i = 1,2,\dots,n-1\) and $F(x) = \mu(-\infty,x]$ then the $p_i$'s are given by 
\begin{equation*}
    p_{i}=\begin{cases}
        F\left(m_{1}\right) & i=1\\
        F\left(m_{i}\right)-F\left(m_{i-1}\right) & i=2,\dots,n-1\\
        1-F\left(m_{n}\right) & i=n,
        \end{cases}
\end{equation*}
where the $x_i$'s solve 
\begin{equation*}
    2F\left(x_{i}\right)=\begin{cases}
    F\left(m_{1}\right) & i=1\\
    F\left(m_{i}\right)+F\left(m_{i-1}\right) & i=2,\dots,n-1\\
    1+F\left(m_{n-1}\right) & i=n.
    \end{cases}
\end{equation*}

\subsection{Asymptotically optimal quantizer for Wasserstein-1}
Finally, we may also compute the asymptotically-optimal representation \cite[Section 7.3]{MR1764176}, \cite[Condition (C')]{MR730904}. We first assume that $ \id \mu(x) = f(x) \id x $ with the following additional assumptions (with $r=1$ in our setting)
\begin{equation}\label{e.asymp-opt-assumptions}
    \begin{split}
    &\int |x| f(x) \id x < \infty,\\
    &\int f(x)^{1/2} \id x < \infty, \\
    \lim_{y \to \infty} 
    &\frac{
        \int_y^\infty |x-y| f(x) \id x
        }
        {
            \left(
            \int_y^\infty f(x)^{1/2} \id x
            \right)
        } = 0, \\
    \lim_{y \to -\infty} 
    &\frac{
        \int_{-\infty}^y |x-y| f(x) \id x
        }
        {
            \left(
            \int_{-\infty}^y  f(x)^{1/2} \id x
            \right)
        } = 0.
    \end{split}
\end{equation}
The first condition is the standing assumption of finite mean. The second condition ensures that the measure $\mu_{1/2}$ defined below is well-defined, and arises naturally since the empirical measure of asymptotically optimal quantization centers converges to $\mu_{1/2}$ (see \cite[Section 7.1]{MR1764176}). The last two conditions are technical tail conditions required for the asymptotic analysis.
Define  
\begin{equation}
\id \mu_{1/2}(x) = \frac{1}{\int \sqrt{f(t)} \id t} \sqrt{f(x)} \id x ,
\end{equation}
and let $F_{1/2}(x) = \mu_{1/2}(-\infty,x]$ denote the CDF.
Then the asymptotically-optimal positions are given by $x_i = F_{1/2}^{-1}((2i-1)/(2n))$ and the $p_i$'s are defined as in the optimal case.
Worth remarking is that these are \emph{sufficient} conditions on $f$ and other such conditions exist \cite{MR1275992,linder1991asymptotically}.

\section{A quantization algorithm for distributions}\label{s.algorithm}
In this section, we introduce the algorithm at the core of this article. The algorithm is formally a function $T : \CM \times \BZ_{\geq 0} \to \CM$ that for each (continuous) $\mu \in \CM$ and $n \geq 0$ outputs a purely discrete probability measure of $2^n$ atoms.
\begin{definition}\label{def.split-function}
Let $\mu$ be a probability measure on the measure space $(\BR, \CB(\BR))$. We say that any function
\begin{equation}
    f : \CM(\BR) \to \BR,
\end{equation}
is a \emph{split function} if it is continuous with respect to the weak topology on $\CM$ and satisfies $f(\mu) \in  [\inf\supp(\mu),\sup \supp(\mu)]$. 
\end{definition}

\begin{remark}
    Before we proceed with describing the algorithm, we present some natural examples of split functions:
\begin{enumerate}
    \item The \emph{median}, \( f(\mu) = \med(\mu) \).
    \item The \emph{mean}, \( f(\mu) = \bar{\mu} \).
    \item The \emph{geometric mean} (assuming \( \mu \) is supported on \( [0,\infty) \)),
    \begin{equation}
        f(\mu) = 
        \exp\left(
            \int_0^\infty \log(t) \, \mathrm{d}\mu(t)
        \right).
    \end{equation}
\end{enumerate}
It is straightforward to verify that these examples satisfy the required assumptions.
\end{remark}

We now turn to defining the algorithm. Let $\Omega_\emptyset = \BR$ and define
\[
\Omega_- = \{x \in \Omega_\emptyset : x \leq f(\mu)\},\ \Omega_+ = \{x \in \Omega_\emptyset : x > f(\mu)\}
\]
as well as (assuming that $\mu(\Omega_\pm)>0$)
\[
\mu_\pm(A) =    \frac{
    \mu(A \cap \Omega_\pm)  }{\mu(\Omega_\pm)},\ A \in \CB(\BR).
\]
The algorithm $T^f$ can then be compactly defined by the following. 
\begin{definition}\label{def.algorithm}    
Let $T^f(\mu,0) = \delta_{f(\mu)}$ and for $n\geq 1$ let 
\begin{equation}\label{e.algorithm-definition}
    \mu^{(n),f} :=T^f(\mu,n) = \mu(\Omega_-)T^f(\mu_-,n-1)+\mu(\Omega_+)T^f(\mu_+,n-1)
\end{equation}
and we let $F^{(n),f}(x) := \mu^{(n),f}(-\infty,x]$ denote the corresponding CDF.  If $\mu(\Omega_\pm) = 0$ the procedure terminates prematurely and sets $\mu^{(n),f}=\delta_{f(\mu_\mp)},\ \forall n\geq 0$.
\end{definition}

Unraveling the recursive definition one obtains the following more explicit description.
\begin{lemma}\label{l.non-recursive}
    Let $n\in \BZ_{\geq 0}$ and let $\alpha \in \{\pm\}^n$.
    Define $\Omega_\pm,\mu_\pm$ as above and for $\alpha \in \{\pm\}^{n}, n\geq 1$ define $\Omega_{\alpha\pm}$ recursively by
    \[
    \Omega_{\alpha-} = \{x \in \Omega_\alpha : x \leq f(\mu_\alpha)\},\ \Omega_{\alpha+} = \{x \in \Omega_\alpha : x > f(\mu_\alpha)\}    
    \]
    and (assuming that $\mu(\Omega_{\alpha})>0$ for all $\alpha \in \{\pm\}^j, j = 1,2,...,n-1$)
    \[
        \mu_{\alpha\pm}(A) =    
        \frac{ \mu(A \cap \Omega_{\alpha \pm})  }{\mu(\Omega_{\alpha \pm})},\ A \in \CB(\BR)   
    .
    \]
    Then
    \begin{equation}\label{def.approx_meas}
    \mu^{(n),f} = T^f(\mu,n)=
    \sum_{\alpha \in \{\pm\}^n} 
        \mu(\Omega_{\alpha}) \delta_{f(\mu_\alpha)},\  n \geq 0,
    \end{equation}
    and
    \begin{equation}\label{def.approx_cdf}
    F^{(n),f}(x) := \sum_{ \substack{ \alpha \in \{\pm\}^n : \\ f(\mu_\alpha) \leq x  } }   \mu(\Omega_{\alpha}).
\end{equation}
\end{lemma}
\begin{proof}
    The proof is a simple inductive argument. Indeed for $n=0$ the two coincide for any $\mu \in \CM_1$. Now assume that for some $n\geq 1$ we have 
    \[
    T^f(\mu,n)
    =
    \sum_{\alpha \in \{\pm\}^n} 
        \mu(\Omega_{\alpha}) \delta_{f(\mu_\alpha)},    
    \]
    for any $\mu \in \CM_1$. 
    Then by Definition~\ref{def.algorithm} 
    \begin{align*}        
    T^f(\mu,n+1) &= \mu(\Omega_-) T^f(\mu_-,n) +\mu(\Omega_+) T^f(\mu_+,n) \\
    &= 
    \mu(\Omega_-)
        \sum_{\alpha \in \{\pm\}^n} 
            \mu_-(\Omega_{\alpha}) \delta_{f( (\mu_-)_{\alpha})}
   +\mu(\Omega_+)\sum_{\alpha \in \{\pm\}^n} 
            \mu_+(\Omega_{\alpha}) \delta_{f((\mu_+)_{\alpha})}.
    \end{align*}
    Now, note that we have by definition 
    \[
    \mu(\Omega_\pm)\mu_\pm(\Omega_{\alpha}) = \mu(\Omega_{\pm\alpha}), 
    \]
    and 
    \[
    f( (\mu_\pm)_{\alpha}) = f( \mu_{\pm\alpha}),
    \]
    which implies the first part of the Lemma: 
    \begin{equation}
         T^f(\mu,n+1) 
         = 
         \sum_{\alpha \in \{\pm\}^{n}}\mu(\Omega_{-\alpha}) \delta_{f( \mu_{-\alpha})}+ \mu(\Omega_{+\alpha}) \delta_{f( \mu_{+\alpha})}
         =\sum_{\alpha \in \{\pm\}^{n+1}}\mu(\Omega_{\alpha}) \delta_{f( \mu_{\alpha})}.
    \end{equation}
    The second claim follows immediately from the definition of $F^{(n)}$. 
\end{proof}
As a remark, we note that $\Omega_\alpha, \alpha \in \{\pm\}^n$ forms a partition of $\Omega.$ Furthermore, when it is clear from context we will drop $f$ as a superscript only writing $\mu^{(n)},\ T(\mu,n),\ F^{(n)}$. Finally, whenever it is possible, we will refrain from using the notation $T(\mu,n)$ and instead we shall use the notation $\mu^{(n)}$ as a shorthand. 

\begin{example}
    In this example we compute the Wasserstein distance between the distribution of a uniform random variable and its approximation using the recursive algorithm, when the split function is given by the median (note that in this case both the median and mean will coincide throughout the algorithm). In other words, let $\mu(\id x) = \indicator_{[0,1]} \id x $ and we aim to compute $W_1(\mu, \mu^{(n)}).$ It should be somewhat clear that in this case we have 
    \[
    \mu(\Omega_\alpha) = \frac{1}{2^n}, \forall \alpha\in\{\pm\}^n,    
    \]
    and if we order the atoms $(x_i^n)_{i=1}^{2^n} = (f(\mu_\alpha))_{\alpha \in \{\pm\}^n}, n \geq 0$ increasingly we have 
    \[
    x_i^n = \frac{2i-1}{2^{n+1}},\ i = 1,2,\dots,2^n.
    \]
    Hence,
    \begin{align*}
    W_1(\mu, \mu^{(n)}) 
    &= 
    \int_0^{\frac{1}{2^{n+1}}} x \id x 
    + \sum_{i=1}^{2^n-1}
        \int_{\frac{2i-1}{2^{n+1}}}^{\frac{2i+1}{2^{n+1}}}
        \left|
            x-\frac{i}{2^n}
        \right| 
        \id x 
    +\int_{1-\frac{1}{2^{n+1}}}^1 1-x \id x \\
    &= 
    2 \int_0^{\frac{1}{2^{n+1}}} x \id x 
    +\sum_{i=1}^{2^n-1} 
        2\int_{\frac{i}{2^n}}^{\frac{i}{2^{n}}+\frac{1}{2^{n+1}}} 
            x -\frac{i}{2^n}
            \id x\\
    &= 
    2 \int_0^{\frac{1}{2^{n+1}}} x \id x 
    +\sum_{i=1}^{2^n-1} 
        2\int_{0}^{\frac{1}{2^{n+1}}} 
            x \id x\\
    &= 
    \frac{1}{2^{2(n+1)}}
    +
    \left(
        2^{n}-1
    \right)
    \frac{1}{2^{2(n+1)}} = \frac{1}{2^{n+2}}.
    \end{align*}
\end{example}

We end this section with a couple of  remarks. 
\begin{remark}\label{r.lower_bound}
\noindent
\begin{itemize}
    \item 
    While this article is mainly focused on deriving an upper bound we remark that one can immediately obtain a lower bound on the Wasserstein distance owing to the Lipschitz-continuity of the mean. 
    Indeed, note that for any two probability measures $\mu,\nu \in \CM_1$
    \[
    |\bar{\mu}-\bar{\nu}| =
    \left| 
    \int t \id \mu -\int t \id \nu
    \right| 
    \leq 
    W_1(\mu,\nu).
    \]
    Hence the Wasserstein distance is always bounded from below by how well one is able to approximate the mean : 
    \[
    |\bar{\mu}-\bar{\mu}^{(n)}| \leq 
    W_1(\mu,\mu^{(n)}).
    \]
    \item 
    A related observation to the item above is that $f(\mu) = \bar{\mu}$ actually renders this lower-bound trivial. Indeed, suppose that $f(\mu)=\bar{\mu}$ and let $X,X^{(n)}$ be distributed as $\mu$ and $\mu^{(n)}$ respectively. Then by the law of total probability 
    \begin{equation}
    \BE [X^{(n)}] = \sum_{\alpha \in \{\pm\}^n} \mu(\Omega_\alpha) \BE[X | X \in \Omega_\alpha] = \BE[X]
    \end{equation}
    which implies that for any $n\geq 0 $ we have 
    \[
    |\bar{\mu}-\bar{\mu}^{(n)}| =0.
    \]
    \item The median-split function can be interpreted as the greedy choice which minimizes the local error at every iteration. This is due to the fact that the median corresponds to the $L^1$-barycenter:
    \[
    \BE[|X-\med(\mu)|] = \inf_{x\in \BR}\BE[|X-x|],
    \]
    while the mean corresponds to the $L^2$-barycenter:
    \[
    \BE[|X-\bar{\mu}|^2] = \inf_{x\in \BR}\BE[|X-x|^2].
    \]
    \item The mean-split function also corresponds to the TTR algorithm, \cite{9756254}.
\end{itemize}
\end{remark}

\section{Analysis of the algorithm}\label{s.analysis}
This section presents the main theoretical results of the article. In Section~\ref{ss.upper_bound}, we establish the central theorem, and in Section~\ref{ss.corollaries}, we provide two concrete examples of split functions that satisfy the key assumption~\eqref{e.split-fun-assump}.
Throughout the remainder of this article we will assume (unless stated otherwise) that
\begin{equation}
    \mu \text{ is a continuous probability measure.}
\end{equation}
Recall that continuous means that the CDF $F_{\mu}(x)$ is continuous. 
\subsection{Preliminary results}
The lemmas which follow are necessary for the development presented in the remainder of this article. The first one is an elementary inequality.
\begin{lemma}\label{l.sqrt-inequality}
    For any $a,b,c,d >0$ we have
    \begin{equation}
        \sqrt{ab}+\sqrt{cd} \leq \sqrt{(a+c)(b+d)}.
    \end{equation}    
\end{lemma}
\begin{proof}
Squaring both sides we get 
\begin{equation}
    ab+2 \sqrt{abcd} +cd \leq (a+c)(b+d) = ab+ad+cb+cd,
\end{equation}
which yields 
\begin{equation}
    2\sqrt{abcd} \leq ad+cb
\end{equation}
which is equivalent to the AM-GM inequality which completes the proof.
\end{proof}

\begin{lemma}\label{l.1point-minimizer} 
Let $\mu \in \CM_1$ then for any $x\in \BR$ we have
\begin{equation}
    W_1(\mu,\delta_x) = \int_\BR |t-x| \id \mu(t).
\end{equation}
Moreover, if $f$ is a split function and letting 
\begin{equation}
    \bar{\mu}_{+} = \BE_\mu[X | X > f(\mu)], \quad
    \bar{\mu}_{-} = \BE_\mu[X | X \leq f(\mu)],
\end{equation} 
then
\begin{equation}\label{e.split-function-1-point-error}
    \begin{split}        
    W_1( \mu,\delta_{\med(\mu)})  
    & = \frac{1}{2} 
    \left(
        \bar{\mu}_+-\bar{\mu}_-
    \right)\\
    W_1( \mu,\delta_{\bar{\mu}})
    &= 2 \mu(\Omega_-) \left(\bar{\mu} - \bar{\mu}_-\right)  
    = 2 \mu(\Omega_+) \left(\bar{\mu}_+ - \bar{\mu}\right).
    \end{split}    
\end{equation}

\end{lemma} 
\begin{proof}
    From the Wasserstein-1 definition, it follows immediately that for any coupling of $\mu$ and $\delta_x$ we must have 
    \begin{equation*}
    \int |t-s| \id  \gamma(s,t) = \int |t-x| \id \mu(t) = \BE_\mu[|X-x|].
    \end{equation*}
    For the second claim, it is effectively trivial to observe that 
    \begin{equation*}
        W_1( \mu,\delta_{\med(\mu)})  
        = \frac{1}{2} 
        \left(
            \BE
                \left[ 
                    X | X > \med(\mu)  
                \right]-
            \BE
                \left[ 
                    X | X \leq \med(\mu)  
                \right] 
        \right)
    \end{equation*}
    since $\mu(\Omega_-) = \mu(\Omega_+) = \frac{1}{2}$ by the definition of the median. For the mean we first observe  
    \begin{equation*}
    W_1( \mu,\delta_{\bar{\mu}}) =    \mu(\Omega_-) (\bar{\mu}-\bar{\mu}_-) + \mu(\Omega_+) (\bar{\mu}_+-\bar{\mu}).
    \end{equation*}
    Now we observe that 
    \begin{equation}\label{e.mean-decomp}    
        \bar{\mu} = \mu(\Omega_-) \bar{\mu}_- + \mu(\Omega_+) \bar{\mu}_+,
    \end{equation}
    which yields the result by repeatedly applying~\eqref{e.mean-decomp} 
    \begin{eqnarray*}
        W_1( \mu,\delta_{\bar{\mu}}) 
        &=& \mu(\Omega_-) (\bar{\mu}-\bar{\mu}_-) + \mu(\Omega_+) (\bar{\mu}_+-\bar{\mu}) \\
        &=& \mu(\Omega_-) \left[\mu(\Omega_-) \bar{\mu}_- + \mu(\Omega_+) \bar{\mu}_+-\bar{\mu}_-\right] \\
        & & \quad +\ \mu(\Omega_+) \left[\bar{\mu}_+-(\mu(\Omega_-) \bar{\mu}_- + \mu(\Omega_+) \bar{\mu}_+)\right]\\
        &=& 2 \mu(\Omega_+)\mu(\Omega_-)(\bar{\mu}_+-\bar{\mu}_-).
    \end{eqnarray*}
    Finally, by either keeping the terms $\bar{\mu}_+$ or $\bar{\mu}_-$ we may also write this expression using~\eqref{e.mean-decomp} as
    \begin{equation*}
    W_1( \mu,\delta_{\bar{\mu}}) = 2 \mu(\Omega_-) (\bar{\mu} - \bar{\mu}_-)  = 2 \mu(\Omega_+) (\bar{\mu}_+ - \bar{\mu}).
    \end{equation*}
\end{proof}

\begin{lemma}\label{l.wasserstein-distance-approximation}
    \begin{equation}
        W_1(\mu,\mu^{(n)}) = \sum_{\alpha \in \{\pm\}^n} \mu(\Omega_\alpha) W_1 ( \mu_\alpha , \delta_{f(\mu_\alpha)})
    \end{equation}
\end{lemma}
\begin{proof}
    The proof is a simple computation using the definition of the Wasserstein-1 metric and Lemma~\ref{l.non-recursive}. Let $F(x)= F_\mu(x)$ for brevity, then
    \begin{align*}
        W_1(\mu,\mu^{(n)}) &= \int_\BR \big|F(x) -F^{(n)}(x)\big| \id x  \overset{\eqref{def.approx_cdf}}{=}\int_\BR \big|F(x) -\sum_{ \substack{ \alpha \in \{\pm\}^n : \\ f(\mu_\alpha) \leq x  } }   \mu(\Omega_{\alpha})\big| \id x  \\ 
    \end{align*}
    First, note that by ordering the atoms increasingly $(f(\mu_\alpha))_{\alpha\in \{\pm\}^n} = (x_j)_{j=1}^{2^n}$ (together with the sets $(\Omega_\alpha)_{\alpha \in \{\pm\}^n} =  ((l_j,u_j])_ {1 \leq j  \leq  2^n -1 },$ where $u_{j-1} = l_j$ and $l_1 = -\infty,\ \Omega_{2^n} = (u_{2^n},\infty)$) we have 
    \begin{eqnarray*}
        \int_\BR \big|F(x) -\sum_{ \substack{ \alpha \in \{\pm\}^n : \\ f(\mu_\alpha) \leq x  } }   \mu(\Omega_{\alpha})\big| \id x
        = \sum_{\alpha \in \{\pm\}^n} \int_{\Omega_\alpha} \big|F(x) -\sum_{ \substack{ \beta \in \{\pm\}^n : \\ f(\mu_\beta) \leq x  } }   \mu(\Omega_{\beta})\big|\id x.
    \end{eqnarray*}
    Now, observe that for any 
    \(
    x \in [x_j,u_j]
    \)
    we have $F^{(n)}(x)  = F(u_j)$ while for  any $x \in (u_{j-1},x_j) $ we have $F^{(n)}(x)  = F(u_{j-1}) = F(l_j)$. Hence, by Fubini's theorem we have 
    \begin{align*}
        &\int_{\Omega_j} |F(x)-F^{(n)}(x)| \id x = \int_{x_j}^{u_j} F(u_j) -F(x) \id x+ \int_{l_j}^{x_j} F(x) - F(u_{j-1}) \id x \\
        &=  \int \indicator\{ x_j \leq x \leq u_j , x < t \leq u_j \} \id \mu(t) \id x  
        +   \int \indicator\{ l_j \leq x \leq x_j , l_{j} < t \leq x \} \id \mu(t) \id x \\
        &= \int_{x_j}^{u_j}  (t-x_j) \id \mu(t)+ \int_{l_j}^{x_j} (x_j -t) \id \mu(t) = \BE[|X- x_j| ; X \in \Omega_j] = \mu(\Omega_j) W_1(\mu_j,\delta_{x_j}).
    \end{align*}
    The argument above extends for $\Omega_1,\ \Omega_{2^n}$ by using the fact that $\mu \in \CM_1$ which implies 
    \begin{equation*}
    W_1(\mu,\mu^{(n)}) = \sum_{j=1}^{2^n}\mu(\Omega_j) W_1(\mu_j,\delta_{x_j}) ,
    \end{equation*}    
    which concludes the argument.
\end{proof}

As a corollary we deduce the following identity:
\begin{corollary}\label{c.recursive-identity}
Let $\mu \in \CM_1$ then for $n \geq 1$ we have
\begin{equation*}
    W_1(\mu,\mu^{(n)}) = \mu(\Omega_-)W_1(\mu_-,\mu_-^{(n-1)}) +\mu(\Omega_+)W_1(\mu_+,\mu_+^{(n-1)}).
\end{equation*}    
\end{corollary}
\begin{proof}
    By Lemma~\ref{l.wasserstein-distance-approximation} we have 
    \begin{align*}
        W_1(\mu,\mu^{(n)}) &= \sum_{\alpha \in \{\pm\}^n} \mu(\Omega_\alpha) W_1 ( \mu_\alpha , \delta_{f(\mu_\alpha)}) \\
        &=\sum_{\alpha \in \{\pm\}^{n-1}} \mu(\Omega_{+\alpha}) W_1 ( \mu_{+\alpha} , \delta_{f(\mu_{+\alpha})}) \\
        &\quad+\ \sum_{\alpha \in \{\pm\}^{n-1}} \mu(\Omega_{-\alpha}) W_1 ( \mu_{-\alpha} , \delta_{f(\mu_{-\alpha})}) \\
        &= \mu(\Omega_+)\sum_{\alpha \in \{\pm\}^{n-1}} \mu_+(\Omega_{+\alpha}) W_1 ( \mu_{+\alpha} , \delta_{f(\mu_{+\alpha})})\\
        &\quad+\ \mu(\Omega_-)\sum_{\alpha \in \{\pm\}^{n-1}} \mu_-(\Omega_{-\alpha}) W_1 ( \mu_{-\alpha} , \delta_{f(\mu_{-\alpha})})\\
        &= \mu(\Omega_-)W_1(\mu_-,\mu_-^{(n-1)}) +\mu(\Omega_+)W_1(\mu_+,\mu_+^{(n-1)}),
    \end{align*}
    where we in the last equality used the fact that $\mu_\pm(\Omega_{\pm\alpha}) = \mu_\pm(\Omega_{\alpha}^{\mu_\pm}).$ The superscript of $\Omega_{\alpha}^{\mu_\pm}$ emphasizes that the recursive algorithm is applied to $\mu_\pm$ (in $n-1$ steps) and not $\mu$.    
\end{proof}

\begin{remark}
    It should be noted that by the convexity of the Wasserstein-1 distance we have for any probability measures $\mu,\nu,\gamma \in \CM_1$ and $\lambda \in [0,1]$  
    \begin{equation*}
    W_1(\mu,\lambda \nu +(1-\lambda)\gamma) 
    \leq 
    \lambda W_1(\mu, \nu ) + (1-\lambda)W_1(\mu, \gamma).
    \end{equation*}
    The above proof then shows that we in fact have equality 
    \begin{eqnarray*}
        W_1( \mu , T(\mu,n)) 
        &\overset{\eqref{e.algorithm-definition}}{=}& 
            W_1
                \left( 
                    \mu , 
                    \mu(\Omega_-)T(\mu_-,n-1)
                    +\mu(\Omega_+)T(\mu_+,n-1)
                \right) \\
        &=& 
            \mu(\Omega_-) W_1(\mu , T(\mu_-,n-1))
            +\mu(\Omega_+) W_1(\mu , T(\mu_+,n-1)).
    \end{eqnarray*}
    In other words, the measure $T(\mu,n)$ is in the set of affine pairs of the convex function $\nu \mapsto W_1(\mu,\nu)$.
\end{remark}

\subsection{A general upper bound for finite mean distributions}\label{ss.upper_bound}
This section contains the main results for finite mean distributions. For simplicity, we state it for measures that are supported on the positive line but the result extends verbatim to $\BR$.
\begin{theorem}\label{thm.unbounded-support-error-bound}
    Let $\mu \in \CM_1$ and assume that $\supp(\mu) = \BR_+$. Assume that $f$ is a split function satisfying 
    \begin{equation}\label{e.split-fun-assump}
    W_1(\nu,\nu^{(n)}) \leq c(f) \frac{b-a}{2^n},\ c(f)>0,
    \end{equation}
    for all probability measures $\nu$ with $\supp(\nu) \subseteq [a,b]$ for some $-\infty<a<b<\infty$.
    
    Let $+{_j} \in \{\pm\}^{j}$ be the element consisting of $j$ consecutive $+$'s (with $+_0 = \emptyset$) and define 
    \[
        \omega_j = f(\mu_{+_j}), j\geq 0,\
        \omega_{-1}=0
    \] (with the convention $\mu_\emptyset = \mu$) and $\Omega_j :=[\omega_{j-1},\ \omega_{j}]$. Then, for any $n\geq 1$, we have
    \begin{equation}\label{e.ttr-error}
        \begin{split}        
            W_1(\mu,\mu^{(n)}) 
            \leq 
            c(f) 
            \sum_{j=0}^{n-1} 
                \frac{
                        (\omega_{j}-\omega_{j-1})\mu(\Omega_{j})
                    }
                    {2^{n-j-1}} 
            + \BE[ |X - \omega_{n}| ; X \geq \omega_{n-1}].
        \end{split}
    \end{equation}
\end{theorem}
\begin{proof}
The proof is an application of Corollary~\ref{c.recursive-identity}.
Let $+{_j}- \in \{\pm\}^{j+1}$ be the element consisting of $j$ consecutive $+$'s and terminating $-$. Then note that by Corollary~\ref{c.recursive-identity} we can deduce by induction the identity
\begin{equation}\label{e.exact_recursive_error}
    \begin{split}
    W_1(\mu,\mu^{(n)}) 
    &= \mu(\Omega_-)W_1(\mu_-,\mu_-^{(n-1)}) 
    +\mu(\Omega_+)W_1(\mu_+,\mu_+^{(n-1)})\\
    &= \sum_{j=0}^{n-1} \mu(\Omega_{+_j-})W_1(\mu_{+_j-},\mu_{+_j-}^{(n-(j+1))}) + \mu(\Omega_{+_n})W_1(\mu_{+_n},\mu_{+_n}^{(0)}).
    \end{split}
\end{equation}
Now, observe that for each $j=0,1,\dots,n-1$ we have 
\[
    \Omega_j := \supp(\mu_{+_j-})=\Omega_{+_j-} = [\omega_{j-1},\ \omega_{j}]
\] 
and hence the assumption applies for each $j=0,1,\dots,n-1$:
\begin{equation*}
    W_1(\mu_{+_j-},\mu_{+_j-}^{(n-(j+1))}) \leq c(f) \frac{\omega_{j}-\omega_{j-1}}{2^{n-(j+1)}}.
\end{equation*} 
Moreover, 
\begin{equation*}
\mu(\Omega_{+_n})W_1(\mu_{+_n},\mu_{+_n}^{(0)}) = \BE[ |X - \omega_{n}| ; X \geq \omega_{n-1}]
\end{equation*}
from which we deduce that 
\begin{equation*}
    \begin{split}
        W_1(\mu,\mu^{(n)}) 
        &= \sum_{j=0}^{n-1} 
            \mu(\Omega_{+_j-})W_1(\mu_{+_j-},\mu_{+_j-}^{(n-(j+1))}) 
            + \BE[ |X - \omega_{n}| ; X \geq \omega_{n-1}]\\
        &\leq c(f) \sum_{j=0}^{n-1} 
        \frac{(\omega_{j}-\omega_{j-1})\mu(\Omega_j)}{2^{n-j-1}}
        + \BE[ |X - \omega_{n}| ; X \geq \omega_{n-1}]
    \end{split}
\end{equation*}
which concludes the proof.
\end{proof}

\begin{remark}
\noindent
    \begin{itemize}
        \item 
        We remark that this theorem provides a general approach. Given a split function, it suffices to prove an upper bound of the Wasserstein-1 distance of order $c(f)(b-a)/2^n$ for distributions $\mu$ with support $\supp(\mu) \subseteq [a,b]$ for some $-\infty < a <b <  \infty$. The above result then implies that this is enough to deduce the general result for arbitrary distributions.
        \item We remark that from Equation~\eqref{e.exact_recursive_error} we can deduce a lower bound, 
        \begin{equation}\label{e.lower_bound}
            W_1(\mu,\mu^{(n)}) \geq 
            \BE 
                \left[
                |X-\omega_{n}| ; X \geq \omega_{n-1}
                \right]
        \end{equation}        
        that in certain scenarios appears to actually capture the correct rate, see Section~\ref{ss.examples} and the proof of Theorem~\ref{thm.poly-tails}.
    \end{itemize}
\end{remark}
As a sanity check we establish that $\mu^{(n)} \to \mu$ with respect to $W_1$.
\begin{corollary}\label{c.consistency}
    Under the assumptions of Theorem~\ref{thm.unbounded-support-error-bound} we have
    \[
    \lim_{n\to \infty} W_1(\mu,\mu^{(n)}) =0.
    \]
\end{corollary}
\begin{proof}
We begin with proving that $\lim_{j\to \infty} \omega_j = \infty.$ To do so we shall show that $(\omega_j)_{j\geq 1}$ is an unbounded strictly increasing sequence.

Now, it should be clear that $(\omega_j)_{j\geq -1}$ forms an increasing sequence, but to be more precise it is in fact strictly increasing. Indeed suppose $\omega_j = \omega_{j-1}$ for some $j\geq 0$ then 
\begin{equation*}    
0=\omega_j-\omega_{j-1} = \BE[X -\omega_{j-1} | X \geq \omega_{j-1}]\geq 0,
\end{equation*}
which implies that $\mu(  X = \omega_{j-1} | X \geq \omega_{j-1}) = 1.$ This is equivalent to 
\begin{equation*}    
\mu( X = \omega_{j-1}) = \mu( X \geq \omega_{j-1}),
\end{equation*}
but since $\mu$ is continuous we have $0=\mu(X = \omega_{j-1})=\mu( X \geq \omega_{j-1})$. Hence we have a contradiction to the assumption that $\supp(\mu) = \BR_+$, thus we conclude $\omega_{j-1} < \omega_j.$ 

Recall that by Definition~\ref{def.split-function}, we have that $f(\nu) \in \supp(\nu), \forall \nu \in \CM$ and $f$ is continuous. Hence, if $\lim_{j\to \infty} \omega_j = a < \infty$,
by recalling the definition of $\omega_j$:
\[
\omega_j = f(\mu_{+_j})
\]
we have by the continuity and the the fact that $f(\mu_{+_j}) \in [\inf \supp(\mu_{+j}),\sup \supp(\mu_{+j})]$ 
\[
\lim_{j\to \infty} f(\mu_{+_j}) = \sup (\supp(\mu)).
\]
Thus, $\lim_{j\to \infty} \omega_j = a < \infty$ contradicts once again the assumption that $\supp(\mu) = \BR_+.$ Now, recall Equation~\eqref{e.ttr-error}
\begin{align*}
    W_1(\mu,\mu^{(n)}) 
            \leq 
            c(f) 
            \sum_{j=0}^{n-1} 
                \frac{
                        (\omega_{j}-\omega_{j-1})\mu(\Omega_{j})
                    }
                    {2^{n-j-1}} 
            + \BE[ |X - \omega_{n}| ; X \geq \omega_{n-1}].
\end{align*}
First we show that $ \BE[ |X - \omega_{n}| ; X \geq \omega_{n-1}] \to 0,\ n\to \infty.$ Note that 
\begin{align*}
    \BE 
    \left[ 
        \left|
            X-\omega_n
        \right|
    \right]
    = \BE[ |X - \omega_{n}| ; X \geq \omega_{n-1}] +\BE[ |X - \omega_{n}| ; X \leq \omega_{n-1}].
\end{align*}
Furthermore, since $\omega_n \to \infty$ there exists some $N >0 $ such that $\omega_n\geq 1$ for all $n \geq N$. Thus
\[
\frac{1}{\omega_n}|X-\omega_n|\indicator\{X \leq \omega_{n-1}\} \leq \frac{1}{\omega_n}|X-\omega_n| \leq |X|+1, \forall n \geq N.
\]
Additionally since $\BE[|X|] < \infty$ we have that $X$ is almost surely finite which implies that
\[
\frac{1}{\omega_n}|X-\omega_n|, \frac{1}{\omega_n}|X-\omega_n|\indicator\{X \leq \omega_{n-1}\} \to 1
\]
almost surely as $n\to \infty$. Hence, by the dominated convergence theorem we have 
\[
    \lim_{n\to \infty} 
    \frac{1}{\omega_n}
    \BE 
    \left[ 
        \left|
            X-\omega_n
        \right|
    \right] 
    = 1
    = \lim_{n\to \infty} 
    \frac{1}{\omega_n}\BE[ |X - \omega_{n}| ; X \leq \omega_{n-1}].
\]
Hence, 
\[
\BE[ |X - \omega_{n}| ; X \geq \omega_{n-1}] \to 0, n\to \infty.
\]
Now we show that the sum tends to $0$. First we note that 
\begin{eqnarray*}
    (\omega_j-\omega_{j-1})\mu(\Omega_j) 
    &=& 
        \BE
            \left[ 
                X-\omega_{j-1} | X \geq \omega_{j-1}
            \right]
        \mu(\Omega_j)\\
    &\leq&
        \BE
            \left[ 
                X-\omega_{j-1} ; X \geq \omega_{j-1}
            \right] =: a_j.
\end{eqnarray*}
Moreover, 
\begin{align*}
    \sum_{j=0}^{n-1} \frac{a_j}{2^{n-j-1}} = \sum_{l=0}^{n-1} \frac{a_{n-1-l}}{2^{l}}
\end{align*}
and note that for each fixed $l$ we have 
\[
\frac{a_{n-1-l}}{2^{l}} \to 0,\ n \to \infty
\]
since $a_j \leq \BE[X;X \geq \omega_j]$ and 
\[
\frac{a_{n-1-l}}{2^{l}}  \leq \frac{\BE[X]}{2^l},
\]
hence we have by the dominated convergence theorem applied to the counting measure on $\BZ_{\geq 0}$ and dominating function $\BE[X]/2^l$ once again that 
\[
\lim_{n\to \infty }
\sum_{j=0}^{n-1} 
                \frac{
                        (\omega_{j}-\omega_{j-1})\mu(\Omega_{j})
                    }
                    {2^{n-j-1}}  = 0,
\]
which completes the proof.
\end{proof}

\subsection{Corollaries}\label{ss.corollaries}
In this section we show that the general bound established in Theorem~\ref{thm.unbounded-support-error-bound} is applicable to the split functions $f(\mu) = \bar{\mu}$ and $f(\mu) = \med(\mu)$. We begin with a technical lemma which is necessary for $f(\mu) = \bar{\mu}$ and can be skipped at a first read-through.
\begin{lemma}\label{l.recursive-upper-bound}
    Let $f(\mu) = \bar{\mu}$, then for any $n\geq 1$ ,
    \begin{equation}\label{e.recursive-upper-bound}    
        \begin{split}            
        &\sum_{\alpha \in \{\pm\}^n } 
            \mu(\Omega_\alpha)  
            \sqrt{
                \left( 
                    \bar{\mu}_{\alpha+} - \bar{\mu}_\alpha 
                \right)
                \left( 
                    \bar{\mu}_\alpha - \bar{\mu}_{\alpha-} 
                \right)
            }\\
        &\leq
        \frac{1}{2}
        \sum_{\alpha \in \{\pm\}^{n-1} } 
            \mu(\Omega_\alpha)  
            \sqrt{
                \left( 
                    \bar{\mu}_{\alpha++} - \bar{\mu}_\alpha 
                \right)
                \left( 
                    \bar{\mu}_\alpha - \bar{\mu}_{\alpha--} 
                \right)
            }.
        \end{split}
    \end{equation}
\end{lemma}
\begin{proof}
The proof is rather straightforward though somewhat cumbersome. First note that by partitioning the sum into strings that end with $+$'s and $-$'s respectively we get
\begin{align*}
    \sum_{\alpha \in \{\pm\}^n } &
        \mu(\Omega_\alpha)  
        \sqrt{
            \left( 
                \bar{\mu}_{\alpha+} - \bar{\mu}_\alpha 
            \right)
            \left( 
                \bar{\mu}_\alpha - \bar{\mu}_{\alpha-} 
            \right)
        }\\
    &=
    \sum_{\alpha \in \{\pm\}^{n-1} } 
        \mu(\Omega_{\alpha+})  
        \sqrt{
            \left( 
                \bar{\mu}_{\alpha++} - \bar{\mu}_{\alpha+} 
            \right)
            \left( 
                \bar{\mu}_{\alpha+} - \bar{\mu}_{\alpha+-} 
            \right)
        }\\
    &\quad +\ \sum_{\alpha \in \{\pm\}^{n-1} }
        \mu(\Omega_{\alpha-})  
        \sqrt{
            \left( 
                \bar{\mu}_{\alpha-+} - \bar{\mu}_{\alpha-} 
            \right)
            \left( 
                \bar{\mu}_{\alpha-} - \bar{\mu}_{\alpha--} 
            \right)
        }\\
        &=
        \sum_{\alpha \in \{\pm\}^{n-1} } 
            \mu(\Omega_{\alpha}) \mu_{\alpha}(\Omega_{\alpha+})  
            \sqrt{
                \left( 
                    \bar{\mu}_{\alpha++} - \bar{\mu}_{\alpha+} 
                \right)
                \left( 
                    \bar{\mu}_{\alpha+} - \bar{\mu}_{\alpha+-} 
                \right)
            }\\
        &\quad + \ \sum_{\alpha \in \{\pm\}^{n-1} }
            \mu(\Omega_{\alpha}) \mu_{{\alpha}}(\Omega_{\alpha-}) 
            \sqrt{
                \left( 
                    \bar{\mu}_{\alpha-+} - \bar{\mu}_{\alpha-} 
                \right)
                \left( 
                    \bar{\mu}_{\alpha-} - \bar{\mu}_{\alpha--} 
                \right)
            }.
\end{align*}
Now note that we clearly have
\begin{equation}\label{e.recursive-mean-bounds}
    \begin{split}
        &\bar{\mu}_{\alpha+} - \bar{\mu}_{\alpha+-} \leq \bar{\mu}_{\alpha+} - \bar{\mu}_{\alpha} \\
        &\bar{\mu}_{\alpha-+} - \bar{\mu}_{\alpha-} \leq \bar{\mu}_{\alpha} - \bar{\mu}_{\alpha-},         
    \end{split}
\end{equation}
which yields 
\begin{equation}\label{e.intermediate-inequality}
    \begin{split}        
    \sum_{\alpha \in \{\pm\}^n } &
        \mu(\Omega_\alpha)  
        \sqrt{
            \left( 
                \bar{\mu}_{\alpha+} - \bar{\mu}_\alpha 
            \right)
            \left( 
                \bar{\mu}_\alpha - \bar{\mu}_{\alpha-} 
            \right)
        } \\
        &\overset{\eqref{e.recursive-mean-bounds}}{\leq} \sum_{\alpha \in \{\pm\}^{n-1} } 
        \mu(\Omega_{\alpha}) \mu_{\alpha}(\Omega_{\alpha+})  
        \sqrt{
            \left( 
                \bar{\mu}_{\alpha++} - \bar{\mu}_{\alpha+} 
            \right)
            \left( 
                \bar{\mu}_{\alpha+} - \bar{\mu}_{\alpha} 
            \right)
        }\\
    &\quad+\ \sum_{\alpha \in \{\pm\}^{n-1} }
        \mu(\Omega_{\alpha}) \mu_{{\alpha}}(\Omega_{\alpha-}) 
        \sqrt{
            \left( 
                \bar{\mu}_{\alpha} - \bar{\mu}_{\alpha-} 
            \right)
            \left( 
                \bar{\mu}_{\alpha-} - \bar{\mu}_{\alpha--} 
            \right)
        }.
    \end{split}
\end{equation}
Now, we observe that by Lemma~\ref{l.1point-minimizer} we have 
\begin{equation*}
W_1(\mu_{\alpha}, \delta_{\bar{\mu}_\alpha})= \mu_{\alpha}(\Omega_{\alpha+})(\bar{\mu}_{\alpha+}-\bar{\mu}_\alpha) = \mu_{\alpha}(\Omega_{\alpha-})(\bar{\mu}_{\alpha}-\bar{\mu}_{\alpha-})
\end{equation*}
which implies that 
\begin{eqnarray}\label{e.sqrt-bounds1}\nonumber
    \mu_{\alpha}(\Omega_{\alpha+})
    \sqrt{
        \left( 
            \bar{\mu}_{\alpha++} - \bar{\mu}_{\alpha+} 
        \right)
        \left( 
            \bar{\mu}_{\alpha+} - \bar{\mu}_{\alpha} 
        \right)
    }
    &=&
    \sqrt{
        \mu_{\alpha}(\Omega_{\alpha+})
        \mu_{\alpha}(\Omega_{\alpha-})
        \left( 
            \bar{\mu}_{\alpha++} - \bar{\mu}_{\alpha+} 
        \right)
        \left( 
            \bar{\mu}_{\alpha} - \bar{\mu}_{\alpha-} 
        \right)
    }\\
    &\leq &
    \frac{1}{2} 
    \sqrt{
        \left( 
            \bar{\mu}_{\alpha++} - \bar{\mu}_{\alpha+} 
        \right)
        \left( 
            \bar{\mu}_{\alpha} - \bar{\mu}_{\alpha-} 
        \right)
    }.
\end{eqnarray}
Analogously,
\begin{eqnarray}
    \label{e.sqrt-bounds2}\nonumber
            \mu_{{\alpha}}(\Omega_{\alpha-}) 
            \sqrt{
                \left( 
                    \bar{\mu}_{\alpha} - \bar{\mu}_{\alpha-} 
                \right)
                \left( 
                    \bar{\mu}_{\alpha-} - \bar{\mu}_{\alpha--} 
                \right)
            }
            &=&
            \sqrt{
                \mu_{{\alpha}}(\Omega_{\alpha-}) 
                \mu_{{\alpha}}(\Omega_{\alpha+}) 
                \left( 
                    \bar{\mu}_{\alpha+} - \bar{\mu}_{\alpha} 
                \right)
                \left( 
                    \bar{\mu}_{\alpha-} - \bar{\mu}_{\alpha--} 
                \right)
            }\\
            &\leq &
            \frac{1}{2}
            \sqrt{
                \left( 
                    \bar{\mu}_{\alpha+} - \bar{\mu}_{\alpha} 
                \right)
                \left( 
                    \bar{\mu}_{\alpha-} - \bar{\mu}_{\alpha--} 
                \right)
            }.
    \end{eqnarray}
Hence, we may bound the right-hand side of~\eqref{e.intermediate-inequality}:
\begin{align*}
    &\sum_{\alpha \in \{\pm\}^{n-1} } 
    \mu(\Omega_{\alpha}) \mu_{\alpha+}(\Omega_{\alpha+})  
    \sqrt{
        \left( 
            \bar{\mu}_{\alpha++} - \bar{\mu}_{\alpha+} 
        \right)
        \left( 
            \bar{\mu}_{\alpha+} - \bar{\mu}_{\alpha} 
        \right)
    }\\
    &+\sum_{\alpha \in \{\pm\}^{n-1} }
    \mu(\Omega_{\alpha}) \mu_{{\alpha-}}(\Omega_{\alpha-}) 
    \sqrt{
        \left( 
            \bar{\mu}_{\alpha} - \bar{\mu}_{\alpha-} 
        \right)
        \left( 
            \bar{\mu}_{\alpha-} - \bar{\mu}_{\alpha--} 
        \right)
    }\\
    &\overset{\eqref{e.sqrt-bounds1},\eqref{e.sqrt-bounds2}}{\leq} 
    \sum_{\alpha \in \{\pm\}^{n-1} } 
    \mu(\Omega_{\alpha})  
    \frac{1}{2} 
        \sqrt{
            \left( 
                \bar{\mu}_{\alpha++} - \bar{\mu}_{\alpha+} 
            \right)
            \left( 
                \bar{\mu}_{\alpha} - \bar{\mu}_{\alpha-} 
            \right)
        }\\
    &+\sum_{\alpha \in \{\pm\}^{n-1} }
    \mu(\Omega_{\alpha})  
    \frac{1}{2}
        \sqrt{
            \left( 
                \bar{\mu}_{\alpha+} - \bar{\mu}_{\alpha} 
            \right)
            \left( 
                \bar{\mu}_{\alpha-} - \bar{\mu}_{\alpha--} 
            \right).
        }
\end{align*}
Finally, by Lemma~\ref{l.sqrt-inequality}, with the following choices
\begin{eqnarray*}
    a &=& \bar{\mu}_{\alpha++} - \bar{\mu}_{\alpha+},\\
    b &=& \bar{\mu}_{\alpha} - \bar{\mu}_{\alpha-},\\
    c &=& \bar{\mu}_{\alpha+} - \bar{\mu}_{\alpha},\\
    d &=& \bar{\mu}_{\alpha-} - \bar{\mu}_{\alpha--},
\end{eqnarray*}
we have
\begin{equation}\label{e.sqrt-inequality-application}
    \begin{split}
    & \sqrt{
            \left( 
                \bar{\mu}_{\alpha++} - \bar{\mu}_{\alpha+} 
            \right)
            \left( 
                \bar{\mu}_{\alpha} - \bar{\mu}_{\alpha-} 
            \right)
        }
    +
    \sqrt{
            \left( 
                \bar{\mu}_{\alpha+} - \bar{\mu}_{\alpha} 
            \right)
            \left( 
                \bar{\mu}_{\alpha-} - \bar{\mu}_{\alpha--} 
            \right)
        }\\
    &\leq
    \sqrt{
        \left(
            (\bar{\mu}_{\alpha++} - \bar{\mu}_{\alpha+}) +(\bar{\mu}_{\alpha+} - \bar{\mu}_{\alpha})
        \right)
        \left(
            (\bar{\mu}_{\alpha} - \bar{\mu}_{\alpha-}) +(\bar{\mu}_{\alpha-} - \bar{\mu}_{\alpha--})
        \right)
    }\\
    &= 
    \sqrt{
        \left(
            \bar{\mu}_{\alpha++}  - \bar{\mu}_{\alpha}
        \right)
        \left(
            \bar{\mu}_{\alpha} - \bar{\mu}_{\alpha--}
        \right)
    },
    \end{split}
\end{equation}
which proves the claim.
\end{proof}

We can now prove our first main result.
\begin{theorem}\label{thm.compact-mean-bound}
    Assume $f(\mu) = \bar{\mu}$ and $\supp(\mu) = [a,b] ,\, -\infty< a < b < \infty$ then 
    \begin{equation}\label{e.compact-mean-bound}
        W_1(\mu,\mu^{(n)})\leq \frac{1}{2} \frac{b-a}{2^{n}}.
    \end{equation}    
\end{theorem}
\begin{proof}
    First we note that by Lemma~\ref{l.1point-minimizer} we have for any $n \geq 1, \alpha \in \{\pm\}^n$ that since 
    \begin{equation*}
    \mu_\alpha(\Omega_{\alpha-}) + \mu_\alpha(\Omega_{\alpha+}) = 1,
    \end{equation*}
    we have 
    \begin{equation*}
    \mu_\alpha(\Omega_{\alpha-}) \mu_\alpha(\Omega_{\alpha+}) \leq \frac{1}{4}.
    \end{equation*}
    Hence, we obtain the bound
    \begin{equation}\label{e.delta-upper-bound}
        \begin{split}
        W_1(\mu_\alpha,\delta_{\bar{\mu}_\alpha}) 
        &=\sqrt{W_1(\mu_\alpha,\delta_{\bar{\mu}_\alpha})^2} \overset{\eqref{e.split-function-1-point-error}}{=} 
        \sqrt{
            4 
            \mu_\alpha(\Omega_{\alpha-})
            \mu_\alpha(\Omega_{\alpha+}) 
            \left( 
                \bar{\mu}_{\alpha+} - \bar{\mu}_\alpha 
            \right)
            \left( 
                \bar{\mu}_\alpha - \bar{\mu}_{\alpha-} 
            \right)
        }\\
        &\leq  \sqrt{
            \left( 
                \bar{\mu}_{\alpha+} - \bar{\mu}_\alpha 
            \right)
            \left( 
                \bar{\mu}_\alpha - \bar{\mu}_{\alpha-} 
            \right)
        }.
        \end{split}
    \end{equation}
    Plugging this inequality in the expression obtained in Lemma~\ref{l.wasserstein-distance-approximation} we get 
    \begin{equation*}
    W_1(\mu,\mu^{(n)}) 
        \leq 
        \sum_{\alpha \in \{\pm\}^n } 
            \mu(\Omega_\alpha)  
            \sqrt{
                \left( 
                    \bar{\mu}_{\alpha+} - \bar{\mu}_\alpha 
                \right)
                \left( 
                    \bar{\mu}_\alpha - \bar{\mu}_{\alpha-} 
                \right)
            }.
    \end{equation*}
    By Lemma~\ref{l.recursive-upper-bound} we get 
    \begin{equation*}
        \begin{split}            
        &\sum_{\alpha \in \{\pm\}^n } 
            \mu(\Omega_\alpha)  
            \sqrt{
                \left( 
                    \bar{\mu}_{\alpha+} - \bar{\mu}_\alpha 
                \right)
                \left( 
                    \bar{\mu}_\alpha - \bar{\mu}_{\alpha-} 
                \right)
            }\\
        &\leq
        \frac{1}{2}
        \sum_{\alpha \in \{\pm\}^{n-1} } 
            \mu(\Omega_\alpha)  
            \sqrt{
                \left( 
                    \bar{\mu}_{\alpha++} - \bar{\mu}_\alpha 
                \right)
                \left( 
                    \bar{\mu}_\alpha - \bar{\mu}_{\alpha--} 
                \right)
            }.
        \end{split}
    \end{equation*}
    By iterating this bound we deduce 
    \begin{equation*}
    W_1( \mu, \mu^{(n)}) 
    \leq 
    \frac{
        \sqrt{
            \left( 
                \bar{\mu}_{+_{n+1}} - \bar{\mu} 
            \right)
            \left( 
                \bar{\mu} - \bar{\mu}_{-_{n+1}} 
            \right)
            }
        }
        {
        2^{n}
        },
    \end{equation*} 
    where $+_{n+1}, -_{n+1}$ denotes the sequences in $\{\pm\}^{n+1}$ consisting of $n+1$ $+$'s and $-$'s respectively. Finally, by the AM-GM inequality we see that 
    \begin{equation*}
        \sqrt{
            \left( 
                \bar{\mu}_{+_{n+1}} - \bar{\mu} 
            \right)
            \left( 
                \bar{\mu} - \bar{\mu}_{-_{n+1}} 
            \right)
            }
            \leq \frac{\bar{\mu}_{+_{n+1}} -\bar{\mu}_{-_{n+1}}}{2} \leq \frac{b-a}{2}. 
    \end{equation*}
\end{proof}
\begin{remark}\label{r.mean-split-remarks}
\noindent
\begin{itemize}
    \item     
    We point out that for certain distributions with unbounded support the proof for the mean-split function already yields a non-trivial bound:
    \begin{equation*}
        W_1(\mu,\mu^{(n)}) \leq \frac{\bar{\mu}_{+_{n+1}} -\bar{\mu}_{-_{n+1}}}{2^{n+1}}.
    \end{equation*}
    In particular, for the Exponential distribution (with unit scale) one can deduce the upper bound 
    \begin{equation*}
        W_1(\mu,\mu^{(n)}) \leq \frac{n+1}{2^{n+1}},
    \end{equation*}
    while for the Pareto distribution with parameter $\alpha \in (1,2]$ this bound does in fact diverge 
    \begin{equation*}
        W_1(\mu,\mu^{(n)}) \leq \frac{\left (\frac{\alpha}{\alpha-1}\right )^{n+1}}{2^{n+1}}.
    \end{equation*}
    We will discuss these two distributions and their error bounds in detail in Section~\ref{ss.examples}.
    \item 
    It is worth noting that, upon inspecting the proofs of Lemmas~\ref{l.1point-minimizer},~\ref{l.wasserstein-distance-approximation},~\ref{l.recursive-upper-bound}, and Theorem~\ref{thm.compact-mean-bound}, one observes that Theorem~\ref{thm.compact-mean-bound} can be extended to discrete distributions. Indeed, if 
    \[
    \mu = \sum_{i=1}^m p_i\delta_{x_i}, \quad x_1 < x_2 < \dots < x_m,
    \]
    then for any \( n = 1, 2, \dots, \lfloor \log_2 m \rfloor \), we can deduce
    \begin{equation} \label{e.discrete_upper_bound}
        W_1(\mu, \mu^{(n)}) \leq \frac{1}{2} \cdot \frac{x_m - x_1}{2^{n}}.
    \end{equation}
    Clearly, this upper bound is not sharp for all values of \( n \). For example, if \( m = 2 \), then by Definition~\ref{def.algorithm} we have
    \[
    T(\mu, 1) = \mu.
    \]
    However, when \( m \) is large and \( x_m - x_1 \) remains bounded, this algorithm can be used to compress the original discrete distribution to a distribution with smaller support without sacrificing much in terms of Wasserstein distance.
    \end{itemize}
\end{remark}

Similarly, for the median-split function we deduce a similar result, though the proof is considerably easier. 
\begin{theorem}\label{thm.compact-median-bound}
    Assume $f(\mu) = \med{(\mu)}$ and $\supp(\mu) = [a,b] , -\infty< a < b < \infty$, then
    \begin{equation*}
        W_1(\mu,\mu^{(n)})\leq \frac{b-a}{2^{n}}.
    \end{equation*}    
\end{theorem}
\begin{proof}
    Note that for any $n\geq 1, \alpha \in \{\pm\}^n $ we have $\mu(\Omega_\alpha) = \frac{1}{2^n}$ which implies  
    \begin{align*}
        W_1(\mu,\mu^{(n)}) = \frac{1}{2^n} \sum_{\alpha \in \{\pm\}^n} \bar{\mu}_{\alpha+}-\bar{\mu}_{\alpha-} \leq \frac{1}{2^n} \sum_{\alpha \in \{\pm\}^n} |\Omega_{\alpha}| = \frac{b-a}{2^n}.
    \end{align*}
\end{proof}

\begin{corollary}\label{c.mean-and-median}
    For $f(\mu)= \bar{\mu}$ and $f(\mu) = \med(\mu)$ we have that the general upper bound of~\eqref{e.ttr-error} holds with 
    \begin{equation*}
        c(\bar{\mu}) = \frac{1}{2},\ c(\med(\mu)) = 1.
    \end{equation*}
\end{corollary}

\subsection{Examples}\label{ss.examples}
\begin{example}
    We first begin with the unit Exponential distribution. Recall that the Exponential distribution (with unit scale parameter) is given by $\mu[x,\infty) = \e^{-x}, x \geq 0$. First recall that by Zador's theorem we have 
    \begin{equation}\label{e.exponential_lower_bound}
        W_1(\mu,\mu^{(n)}) \geq \frac{1+o(1)}{4 \cdot 2^n} \left( \int_0^\infty \e^{-x/2} \right)^2 = \frac{1+o(1)}{2^n}.
    \end{equation}
    We begin with $f(\mu)= \bar{\mu}$. The first thing we need to do is to compute the sequence $\omega_j$. By the memorylessness we get 
    \begin{equation*}
    \omega_{j}= \BE [ X | X \geq \omega_{j-1}] = \omega_{j-1} +1 \Longrightarrow \omega_j = j+1,
    \end{equation*}
    which implies that by Theorem~\ref{thm.unbounded-support-error-bound} together with Corollary~\ref{c.mean-and-median}
    \begin{equation*}
        \begin{split}
        W_1(\mu,\mu^{(n)})
        &\leq 
        \frac{1}{2}
        \sum_{j=0}^{n-1} 
                \frac
                    { (\omega_{j}-\omega_{j-1})\mu(\Omega_{j}) }
                    { 2^{n-j-1} } 
                + \BE[ |X - \omega_{n}| ; X \geq \omega_{n-1}]\\
        &=
        \sum_{j=0}^{n-1} 
        \frac
            { \e^{-\omega_{j-1}}-\e^{-\omega_ {j}}  }
            { 2^{n-j} }
            +\BE[|X-(n+1)|;X\geq n].
        \end{split}
    \end{equation*}    
    Now, 
    \begin{equation}\label{e.exponential_upper_bound_1}
    \begin{split}        
        \sum_{j=0}^{n-1} 
            \frac{ \e^{-\omega_{j-1}}-\e^{-\omega_ {j}} }
                 { 2^{n-j} } 
        &= 
        \frac{1}
             {2^n}
        \sum_{j=0}^{n-1} 2^j
                    \e^{-\omega_{j-1}}
                    \left(
                        1-\e^{-(\omega_ {j}-\omega_{j-1})}
                    \right)
        = 
        \frac{\left( 1-\e^{-1} \right) }
             {  2^n }
        \sum_{j=0}^{n-1}  
            \left(
                \frac{2}{\e}
            \right)^j
        \\
        &= 
        \frac{  
            \left( 
                1-\e^{-1}
                \right) }
            {  2^n }
        \frac{ \e }
             { \e-2 }
            \left(
                1-(2/\e)^n
            \right)
            = \frac{1}{2^n}
        \frac{\e-1}{\e-2}
        \left[
                1-\left(\frac{2}{\e}\right)^n
        \right].
    \end{split}
    \end{equation}
    All that remains is to compute the last term:
    \begin{equation}\label{e.exponential_upper_bound_2}
        \begin{split}
        \BE[|X-n|;X\geq n-1] 
        & = \BE[n-X ; n -1\leq X \leq n]+\BE[X-n; X \geq n] \\
        &= \e^{-n}+\e^{-n} = 2 \e^{-n},
        \end{split}
    \end{equation}
    which implies that
    \begin{equation}\label{e.exponential_upper_bound}
        W_1(\mu,\mu^{(n)}) 
        \overset{\eqref{e.exponential_upper_bound_1},\eqref{e.exponential_upper_bound_2}}{\leq}
        \frac{1}{2^n}
        \frac{\e-1}{\e-2}
        \left[
                1-\left(\frac{2}{\e}\right)^n
        \right]
        +2 \e^{-n}.
    \end{equation}
    Hence, we see that since $1/\e \leq 1/2$, we can deduce that 
    \begin{equation*}
        1 
        \overset{\eqref{e.exponential_lower_bound}}{\leq} 
        \liminf_{n\to \infty} 2^n  W_1(\mu,\mu^{(n)})  
        \leq
        \limsup_{n\to \infty}2^n 
        W_1(\mu,\mu^{(n)}) 
        \overset{\eqref{e.exponential_upper_bound}}{\leq}
        \frac{\e-1}{\e-2} \approx 2.39 .
    \end{equation*}
    For $f(\mu) = \med(\mu)$ we instead have 
    \begin{equation*}
    \omega_j = 
    -\log
        \left(
            1-\frac{2^j-1}{2^j}
        \right)
    = 
    j 
    \log
        \left(
            2
        \right)
    \end{equation*}
    and 
    \begin{equation*}
        (\omega_j-\omega_{j-1})\mu[\omega_{j-1},\ \omega_{j}] = \frac{{\log(2)}}{2^{j}}.
    \end{equation*}
    Thus we obtain a worse upper bound than in the mean-split case, namely:
    \begin{align*}        
        W_1(\mu,\mu^{(n)}) &\leq \sum_{j=0}^{n-1} \frac{\log(2)}{2^j} \frac{1}{2^{n-j-1}}+\BE[|X-\log(2)(n+1)|; X \geq \log(2)n] \\
        &= \frac{2\log(2)n}{2^n} + \frac{ \log(4)-1}{2^{n+1}} +\frac{1}{2^{n+1}}
    \end{align*}
\end{example}

\begin{example}\label{ex.pareto}
    For the Pareto distribution $\mu[x,\infty) = x^{-\alpha}, x\geq 1 ,\ \alpha >1$ we instead have 
    \begin{equation*}
    \omega_j = \frac{\alpha}{\alpha-1}\omega_{j-1}
    \end{equation*}
    which clearly implies that $\omega_j= \left(\frac{\alpha}{\alpha-1}\right)^{j+1}$. For $\alpha\ne2$, a lengthy computation yields that
    \begin{eqnarray*}
            W_{1}
            \left(
                \mu,
                \mu^{\left(n\right)}
            \right)
            &\le&
                \frac{
                     \left(1-\frac{1}{\alpha}\right)^{-\alpha}-1
                     }
                     {
                        2\alpha
                        -\left(\alpha-1\right)\left(1-\frac{1}{\alpha}\right)^{-\alpha}
                     }
                \left[
                    2^{-n}+\left(1-\frac{1}{\alpha}\right)^{n\left(\alpha-1\right)}
                \right]\\
        &&\quad+\ \frac{2}{\alpha-1}\left(1-\frac{1}{\alpha}\right)^{\left(n+1\right)\left(\alpha-1\right)}.
    \end{eqnarray*}
    We can see from this expression that we get two different behaviors. For $\alpha>2$,  $(1-1/\alpha)^{\alpha-1}<2$ and so the
    term $2^{-n}$ dominates. Hence,
    \begin{equation*}
        \limsup_{n\rightarrow\infty}
        2^{n}W_{1}
        \left(\mu,\mu^{\left(n\right)}\right)
        \le
        \frac{\left(1-\frac{1}{\alpha}\right)^{-\alpha}-1}
            {2\alpha-\left(\alpha-1\right)\left(1-\frac{1}{\alpha}\right)^{-\alpha}}.
    \end{equation*}
    For $1<\alpha<2$, we have on the other hand $(1-1/\alpha)^{\alpha-1}>2$ and hence the term $ \left(1-\frac{1}{\alpha}\right)^{n(1-\alpha)}$ dominates:
        \begin{equation*}
        \limsup_{n\rightarrow\infty}\left(1-\frac{1}{\alpha}\right)^{n(1-\alpha)}W_{1}\left(\mu,\mu^{\left(n\right)}\right)\le\frac{\left(1-\frac{1}{\alpha}\right)^{-\alpha}-1}{2\alpha-\left(\alpha-1\right)\left(1-\frac{1}{\alpha}\right)^{-\alpha}}+\frac{2\left(1-\frac{1}{\alpha}\right)^{\alpha-1}}{\alpha-1}.
    \end{equation*}
    Interestingly for $\alpha=2$, we instead get
    \begin{equation*}
    W_{1}\left(\mu,\mu^{\left(n\right)}\right) = \CO\left(n 2^{-n} \right).
    \end{equation*}
    Moreover, it is clear that by Equation~\eqref{e.ttr-error} we have
    \begin{equation*}
    W_1(\mu,\mu^{(n)})\geq \BE
    \left[
        \left|
            X- 
            \left(
                \frac{\alpha}{\alpha-1}
            \right)^{n+1}
        \right|;             
        X\geq  
        \left(
            \frac{\alpha}{\alpha-1}
        \right)^{n}                        
    \right],
    \end{equation*} 
    from which we can conclude that for $\alpha \in (1,2)$ we get 
    $
    W_1(\mu,\mu^{(n)}) 
    = 
    \Theta
    \left(
        \left(
            1-\frac{1}{\alpha}
        \right)^{(\alpha-1)n}
    \right)
    $
    while for $\alpha>2$ we have
    \begin{equation*}
    W_1(\mu,\mu^{(n)}) 
    = 
    \Theta
    \left(
        2^{-n}
    \right).
    \end{equation*}
    More compactly written we have for $\alpha \neq 2$:
    \begin{equation}
    W_1(\mu,\mu^{(n)}) 
    = 
    \Theta
    \left(
        \left(
        \frac{1}{2} 
        \vee 
        \left(
            1-\frac{1}{\alpha}
        \right)^{ \alpha-1}
        \right)^n
    \right).
    \end{equation}
    For the median-split function we instead have 
    \begin{equation*}
    \omega_j = 
    \left(
        1-\frac{2^{j+1}-1}{2^{j+1}}
    \right)^{-1/\alpha} 
    =
    2^{(j+1)/\alpha}
    \end{equation*}
    which now instead yields 
    \begin{equation*}
        (\omega_j-\omega_{j-1})
        \mu[\omega_{j-1},\omega_{j}] 
        \leq 
        c(\alpha) 
        2^{j/\alpha}
        \left(
            \frac{1}{2^{j-1}}
            -
            \frac{1}{2^{j}}
        \right)
        = 
        \frac{
            c(\alpha)
        }{
            2^{(1-\frac{1}{\alpha})j }
        }.
    \end{equation*}
    Similarly, 
    \begin{equation*}
        \BE 
        \left[
            X-\omega_{n}
        ; X\geq \omega_{n-1}
        \right] 
        = \frac{1}{\alpha-1} 2^{j(1-\alpha)/\alpha} = \frac{1}{(\alpha-1)2^{(1-1/\alpha)n}}
    \end{equation*}
    hence the last term is bounded from above by $c(\alpha)  2^{-n(1-1/\alpha)}$. Clearly for any $\alpha>1$ we have $2^{- n} \leq 2^{-(1-\frac{1}{\alpha})n}$ which instead yields a different bound 
    \begin{equation*}
    W_1(\mu,\mu^{(n)}) = \Theta(2^{-(1-\frac{1}{\alpha})n} ), \forall \alpha >1.    
    \end{equation*}
    Moreover, it is an elementary but tedious exercise in analysis to see that 
    \( \left(1-\frac{1}{\alpha}\right)^{\alpha-1} < 2^{-(\alpha-1)}\) which implies that the mean-split function performs better in the entire range $\alpha>1$.

\end{example}

\subsection{Error bounds for distributions with polynomially-decaying tails}
In this section we prove Theorem~\ref{thm.poly-tails}. The proof of this result is essentially a generalization of Example~\ref{ex.pareto}, but the conclusion is weaker due to the presence of $o(1)$ terms appearing in the proof.
We further remark that the result likely holds for regularly-varying distributions but to avoid additional technical details we have chosen to restrict our attention to a smaller class of distributions.
\begin{proof}[Proof of Theorem~\ref{thm.poly-tails}]
    The proof will be performed in two steps. We will prove that 
    \begin{equation}\label{e.poly-tail-upper-lower}
        \begin{split}
            &\limsup_{n \to \infty} \frac{\log(W_1(\mu,\mu^{(n)}))}{n} \leq 
            \log 
            \left(
                \left(
                    1-\frac{1}{\alpha}
                \right)^{(\alpha-1)}
                \vee
                \frac{1}{2}
            \right),
                \\
            &\liminf_{n \to \infty} \frac{\log(W_1(\mu,\mu^{(n)}))}{n} \geq 
            \log 
            \left(
                \left(
                    1-\frac{1}{\alpha}
                \right)^{(\alpha-1)}
                \vee
                \frac{1}{2}
            \right).
            \end{split}
    \end{equation}
    To achieve this, all we need to do is to estimate the difference 
    \begin{equation*}
    \omega_j-\omega_{j-1} = 
    \BE  
    \left[
        X-\omega_{j-1} | X \geq \omega_{j-1}
    \right].
    \end{equation*}
    Thus let $x$ be large, then we have
    \begin{align*}
        \BE  
        \left[
            X-x | X \geq x
        \right] \sim 
        \BE  
        \left[
            X-x ; X \geq x
        \right] \frac{x^{\alpha}}{c} 
        \sim  
        \frac{c x^{1-\alpha}}{\alpha-1}
        \frac{x^{\alpha}}{c}
        = \frac{x}{\alpha-1}.
    \end{align*}
    Hence, for any $\epsilon>0$ there exists some $N\geq 1$ such that 
    \begin{equation*}
        \frac{1-\epsilon}{\alpha-1} \omega_{j-1}
        \leq
        \omega_{j}-\omega_{j-1} 
        \leq 
        \frac{1+\epsilon}{\alpha-1} \omega_{j-1},    
        \forall j\geq N.
    \end{equation*} 
    Therefore, for all $j \geq N$ we get 
    \begin{equation*}
    \omega_N 
    \left(
        \frac{\alpha-\epsilon}{\alpha-1}
    \right)^{j-N}
    \leq 
    \omega_j 
    \leq 
    \omega_N 
    \left(
        \frac{\alpha+\epsilon}{\alpha-1}
    \right)^{j-N},
    \end{equation*}
    which implies that for any $j\geq N$ we have
    \begin{align*}
        &(\omega_{j}-\omega_{j-1}) 
        \mu[\omega_{j-1},\omega_{j}]
        \leq     
        \omega_N 
        \left(
            \frac{\alpha+\epsilon}{\alpha-1}
        \right)^{j-N-1}
        c 
        \omega_N^{-\alpha}
        \left(
            \frac{\alpha-1}{\alpha-\epsilon}
        \right)^{(j-N-1)\alpha}\\
        & = 
        c(N,\alpha) 
        \left(
            \frac{1+\epsilon/\alpha}{(1-\epsilon/\alpha)^\alpha}
        \right)^{j-N-1}
        \left(
            1-\frac{1}{\alpha}
        \right)^{(j-N-1)(\alpha-1)}.
    \end{align*}
    Hence, for $n\geq N$ we have 
    \begin{equation*}
    W_1(\mu,\mu^{(n)}) 
    = 
    \CO 
    \left(
        \left[
        \left(
            \frac{1+\epsilon/\alpha}{(1-\epsilon/\alpha)^\alpha}
        \right)
        \left(
            1-\frac{1}{\alpha}
        \right)^{(\alpha-1)}
        \vee
        \frac{1}{2}
        \right]^n
    \right),
    \end{equation*}
    which implies that
    \begin{equation}\label{e.poly-tail-limsup}
        \limsup_{n \to \infty} \frac{\log(W_1(\mu,\mu^{(n)}))}{n}  \leq 
        \log 
        \left(
            \left(
                \frac{1+\epsilon/\alpha}{(1-\epsilon/\alpha)^\alpha}
            \right)
            \left(
                1-\frac{1}{\alpha}
            \right)^{(\alpha-1)}
            \vee
            \frac{1}{2}
        \right),
        \end{equation}
    but since $\epsilon >0$ was arbitrary the first claim of~\eqref{e.poly-tail-upper-lower} is concluded.
    For the corresponding lower bound, simply observe that by Equation~\eqref{e.ttr-error} we have
    \begin{equation*}
    W_1(\mu,\mu^{(n)}) 
    \geq 
    \BE 
    \left[ 
        \left|    
            X -\omega_{n}
        \right|; 
        X \geq \omega_{n-1} 
    \right]
    \geq 
    \BE 
    \left[ 
        X -\omega_{n};
        X \geq \omega_{n} 
    \right].
    \end{equation*}
    Moreover, by the assumption on $\mu$ we have 
    \begin{equation*}
        \BE 
        \left[ 
            X -\omega_{n};
            X \geq \omega_{n} 
        \right]\sim \omega_{n}^{1-\alpha} 
        \geq c(N,\alpha) 
        \left(
            \frac{\alpha-1}{\alpha+\epsilon}
        \right)^{(\alpha-1)n},
    \end{equation*}
    which implies that 
    \begin{equation}\label{e.poly-tails-liminf}
        \liminf_{n \to \infty} \frac{\log(W_1(\mu,\mu^{(n)}))}{n} 
        \geq 
        \log 
        \left(
            \left(\frac{\alpha-1}{\alpha+\epsilon}\right)^{(\alpha-1)}
         \vee
            \frac{1}{2}
        \right).
    \end{equation}
    Hence, since $\epsilon>0$ was arbitrary we deduce the lower bound of~\eqref{e.poly-tail-upper-lower}. Finally as the right-hand sides of~\eqref{e.poly-tail-limsup} and~\eqref{e.poly-tails-liminf} coincide we deduce the result:
    \begin{equation*}
    \lim_{n\to \infty} \frac{\log(W_1(\mu,\mu^{(n)}))}{n} 
    =
    \log 
    \left(
        \left(
            1-\frac{1}{\alpha}
        \right)^{(\alpha-1)}
        \vee
        \frac{1}{2}
    \right).
    \end{equation*}
\end{proof}
\begin{remark}\label{r.polytail-conclusions}
\noindent
    \begin{itemize}
        \item
            We remark that if $\alpha \geq 2$ then $\mu \in \CM_2$ and
            \begin{equation*}
                W_1(\mu,\mu^{(n)}) = \Theta(2^{-n(1+o(1))}),
            \end{equation*}
            which might indicate that a result of the form
            \begin{equation*}
            \mu \in \CM_r, \ r > 2 \Longrightarrow  \lim_{n\to \infty} \frac{\log(W_1(\mu,\mu^{(n)}))}{n}   = -\log(2)
            \end{equation*}
            could be true. We note that this conjecture is supported by all examples considered in this article. 
            Indeed, for the Exponential distribution (Section~\ref{ss.examples}) and 
            the Pareto distribution with $\alpha > 2$ (Example~\ref{ex.pareto}), 
            the rate $\Theta(2^{-n})$ is achieved.
            Moreover, at the boundary $\alpha = 2$, the rate $\CO(n 2^{-n})$ 
            still yields $\limsup_{n\to\infty} \log(W_1(\mu,\mu^{(n)}))/n = -\log 2$.
        \item
           We note that for the class of distributions considered in Theorem~\ref{thm.poly-tails}, the analogous (but stronger) result holds for the asymptotically optimal quantizers \cite[see the second paragraph at the top of p.~677]{MR730904}. In other words, this class satisfies the conditions given by Equation~\eqref{e.asymp-opt-assumptions}.
        \item
            Another natural family of distributions to study is those with finite moment generating function. In this case it would be reasonable to expect that $W_1(\mu,\mu^{(n)}) = \Theta(2^{-n})$.
    \end{itemize}
\end{remark}

\section{Numerical experiments}\label{s.numerical}
This section contains several numerical experiments relating to the Wasserstein-1 error of various discrete representations applied to various distributions. We evaluate the introduced algorithm and compare it with existing representations such as the optimal and the asymptotically-optimal representations outlined in Section~\ref{s.not-def}. We then investigate how the error changes as we perform arithmetic operations between distributions. In this context, an arithmetic operation between distributions \( \mu \) and \( \nu \) is defined as the binary operation induced by performing the corresponding arithmetic operation on independent random variables \( X \sim \mu \) and \( Y \sim \nu \). For example, \( X + Y \) corresponds to the \emph{additive convolution} $\mu*\nu$, while \( X \cdot Y \) corresponds to the \emph{multiplicative convolution} $\mu \star \nu$.
 
Additionally, we explicitly outline how arithmetic operations can be carried out between an arbitrary number of distributions without suffering from the curse of dimensionality.
\subsection{Closed-form distributions}
We begin with comparing the Wasserstein-1 errors for different distributions and different representation sizes where we know the distribution functions for the desired distribution.

\begin{table}[H]
    \centering
    \caption{Distributions used in the numerical experiments.}
    \label{tab:distributions}
    \begin{tabular}{ll}
        \toprule
        Distribution & Definition \\
        \midrule
        Gaussian & $\mathcal{N}(0, 1)$ \\
        Exponential & $\operatorname{Exp}(1)$ \\
        Pareto & $\operatorname{Pareto}(\alpha = 2, x_m = 1)$ \\
        Heavy-tailed & $\mathbb{P}(X \geq x) = \frac{\e}{x \log^2 x}$, \quad $x \geq \e$ \\
        Bimodal & $f(x) = \frac{2 \e^{-(x-2)^2} + \e^{-(x+2)^2}}{3\sqrt{\pi}}$ \\
        \bottomrule
    \end{tabular}
\end{table}
The distributions in Table~\ref{tab:distributions} were chosen to cover a variety of different cases. In particular, we wanted to test the different representations on distributions with different tail behaviors, symmetry and modality.

\begin{figure}[H]
    \centering
    \begin{subfigure}{0.49\textwidth}
        \centering
        \includegraphics[width=\linewidth]{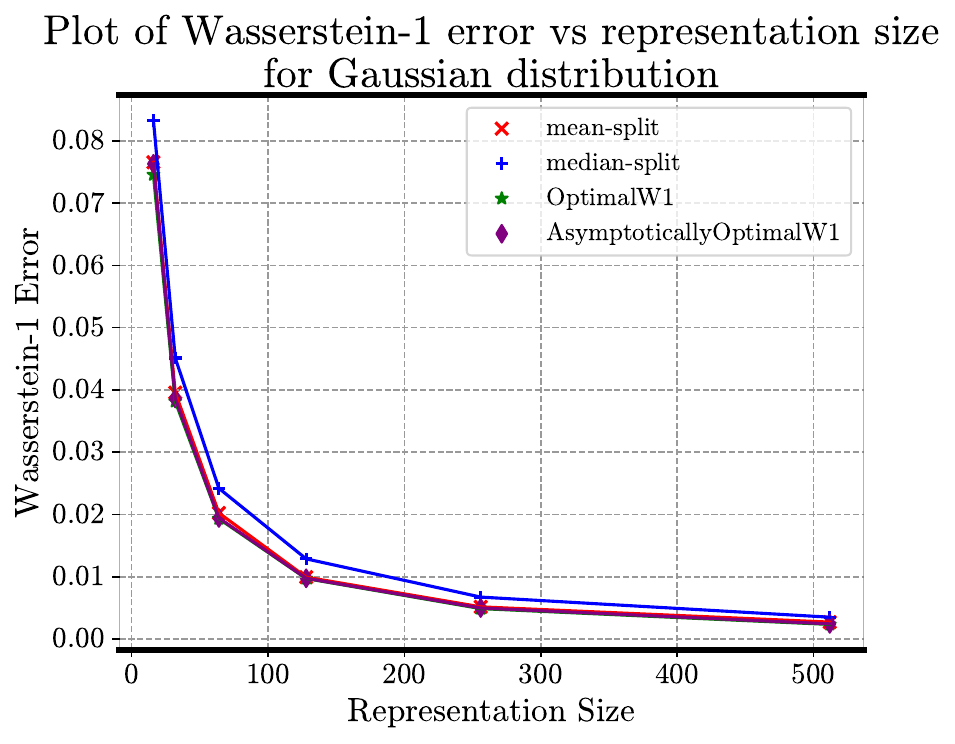}
    \end{subfigure} 
    \begin{subfigure}{0.49\textwidth}
        \centering
        \includegraphics[width=\linewidth]{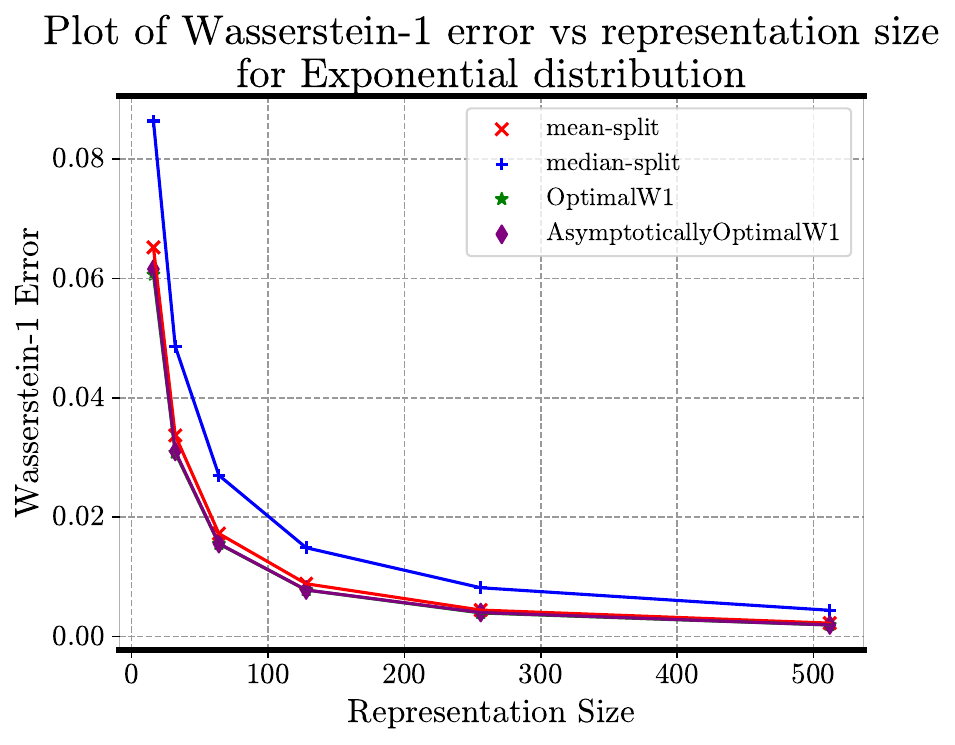}
    \end{subfigure}
    \begin{subfigure}{0.49\textwidth}
        \centering
        \includegraphics[width=\linewidth]{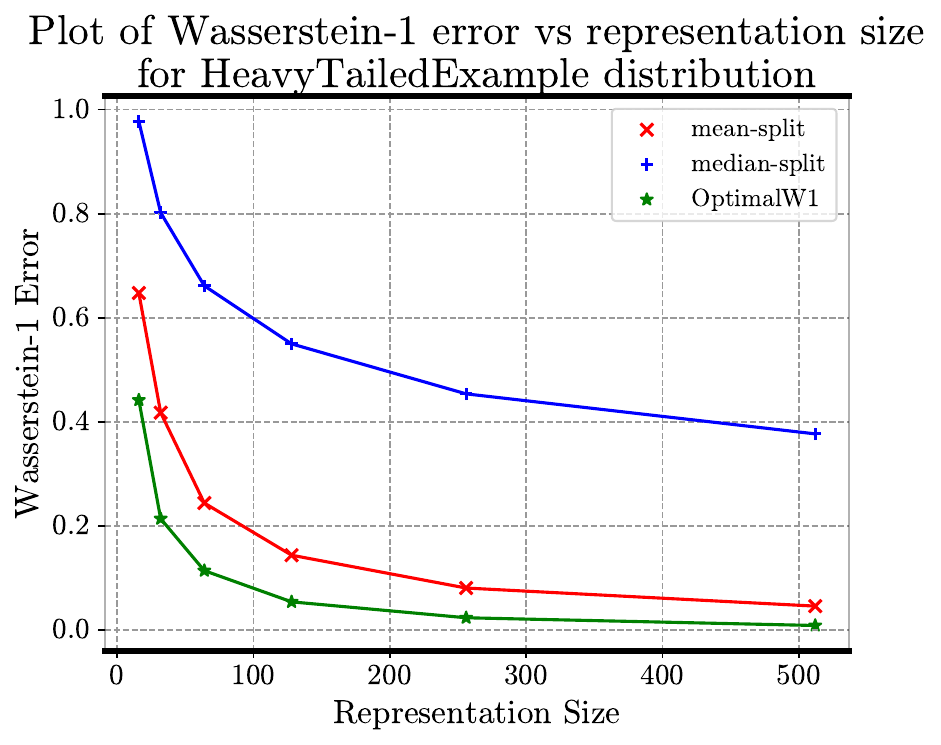}
    \end{subfigure} 
    \begin{subfigure}{0.49\textwidth}
        \centering
        \includegraphics[width=\linewidth]{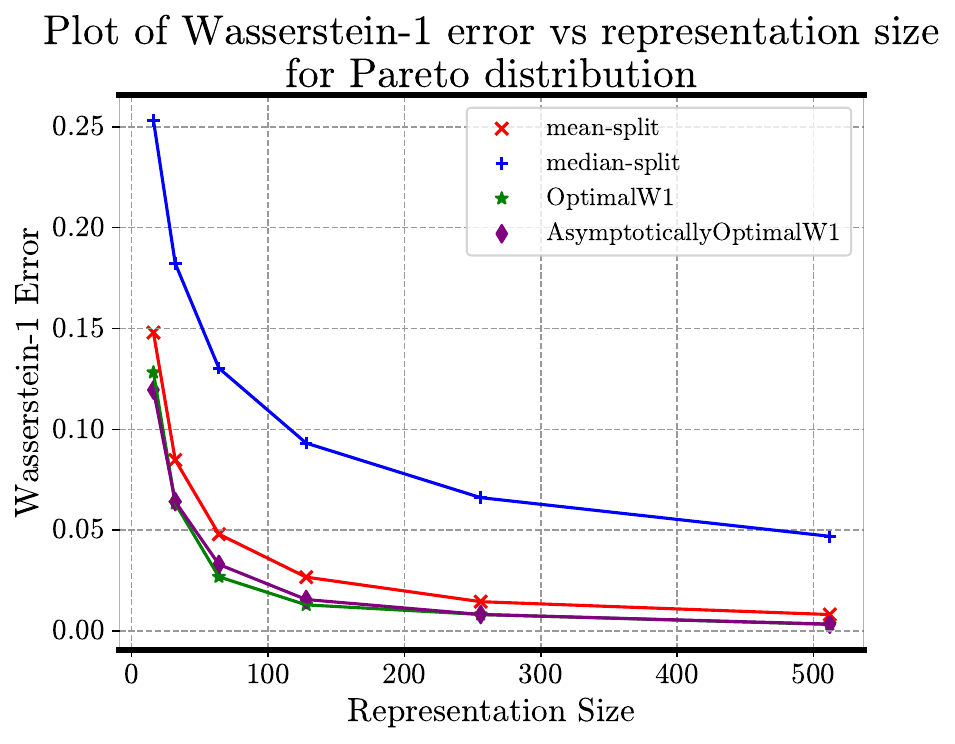}
    \end{subfigure} 
    \begin{subfigure}{0.49\textwidth}
        \centering
        \includegraphics[width=\linewidth]{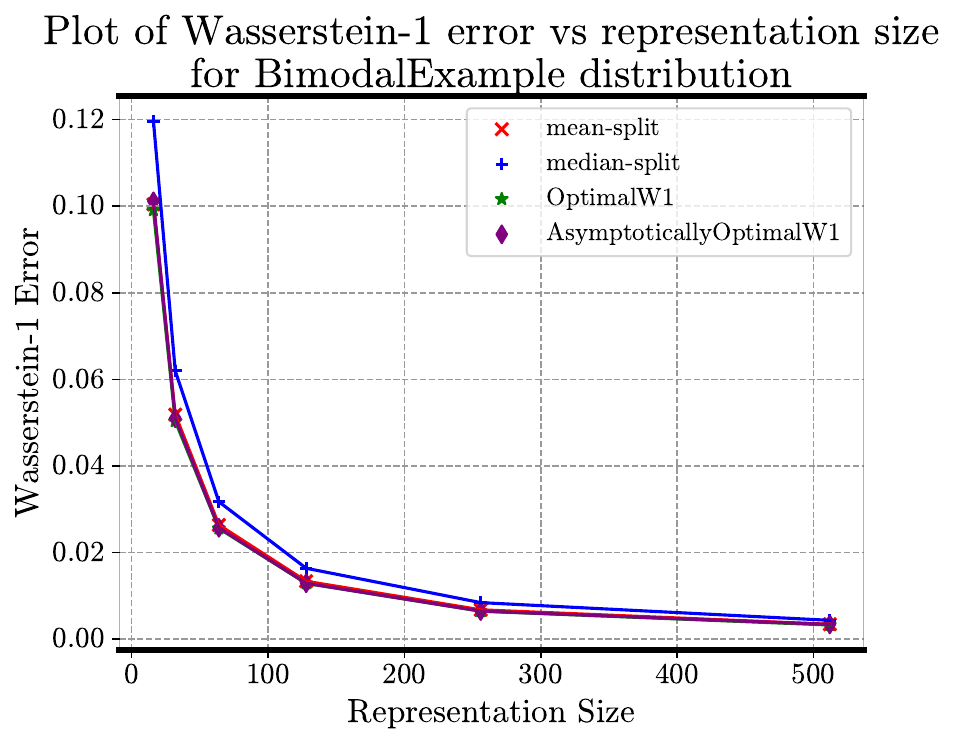}
    \end{subfigure} 
\caption{Wasserstein-1 error vs representation size for mean-split, median-split, optimal and asymptotically optimal representations for the distributions in Table~\ref{tab:distributions}.}    \label{fig:wass_vs_repsize}
\end{figure}

\begin{figure}[H]
    \centering
    \begin{subfigure}{0.49\textwidth}
        \centering
        \includegraphics[width=\linewidth]{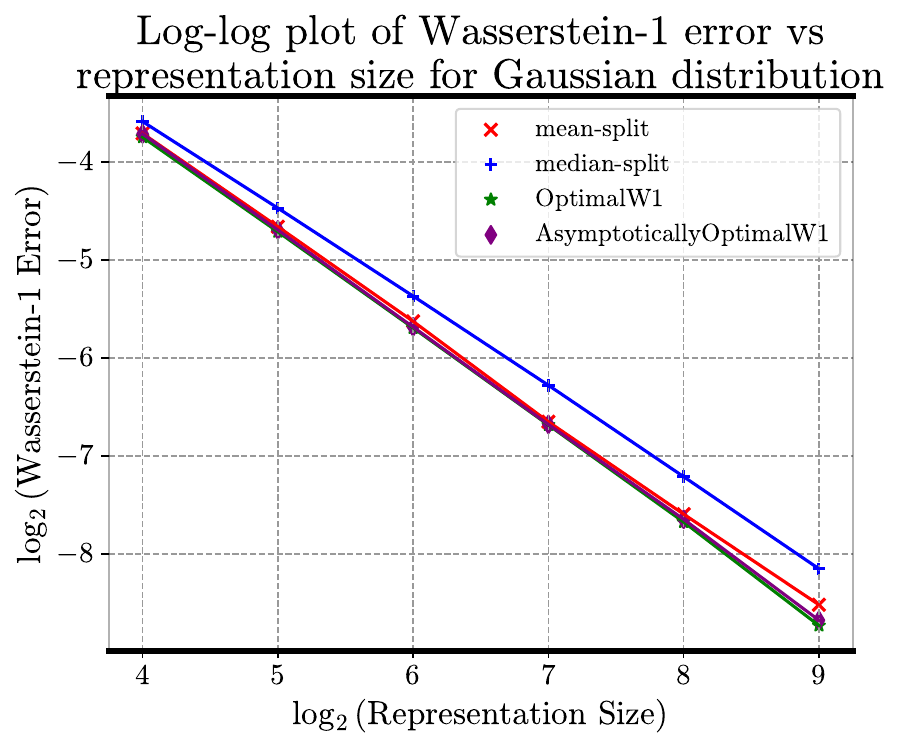}
    \end{subfigure} 
    \begin{subfigure}{0.49\textwidth}
        \centering
        \includegraphics[width=\linewidth]{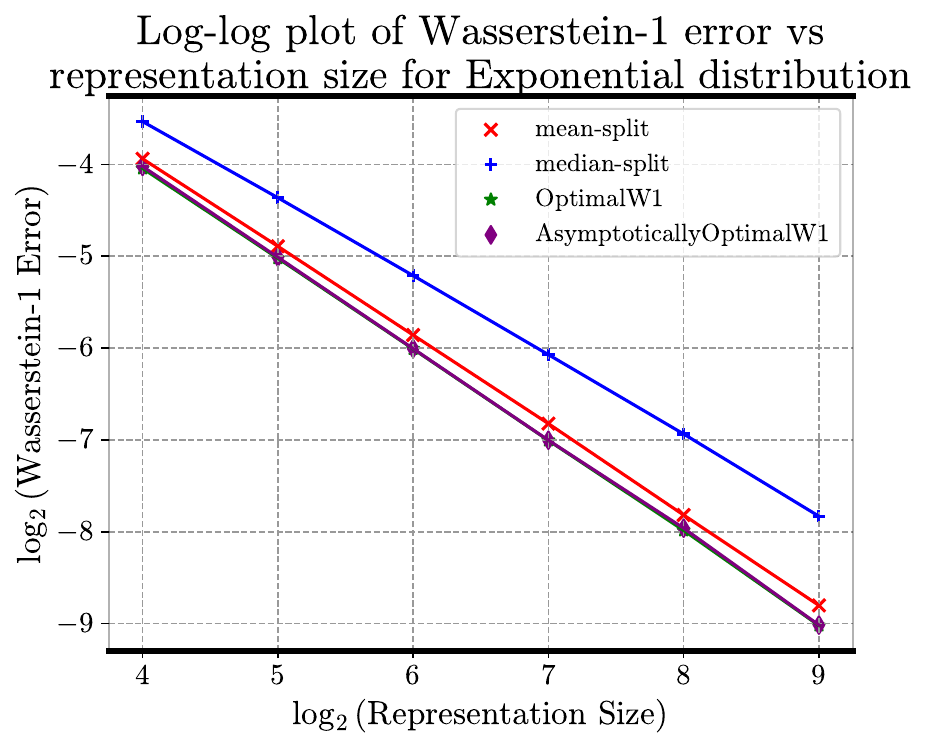}
    \end{subfigure}
    \begin{subfigure}{0.49\textwidth}
        \centering
        \includegraphics[width=\linewidth]{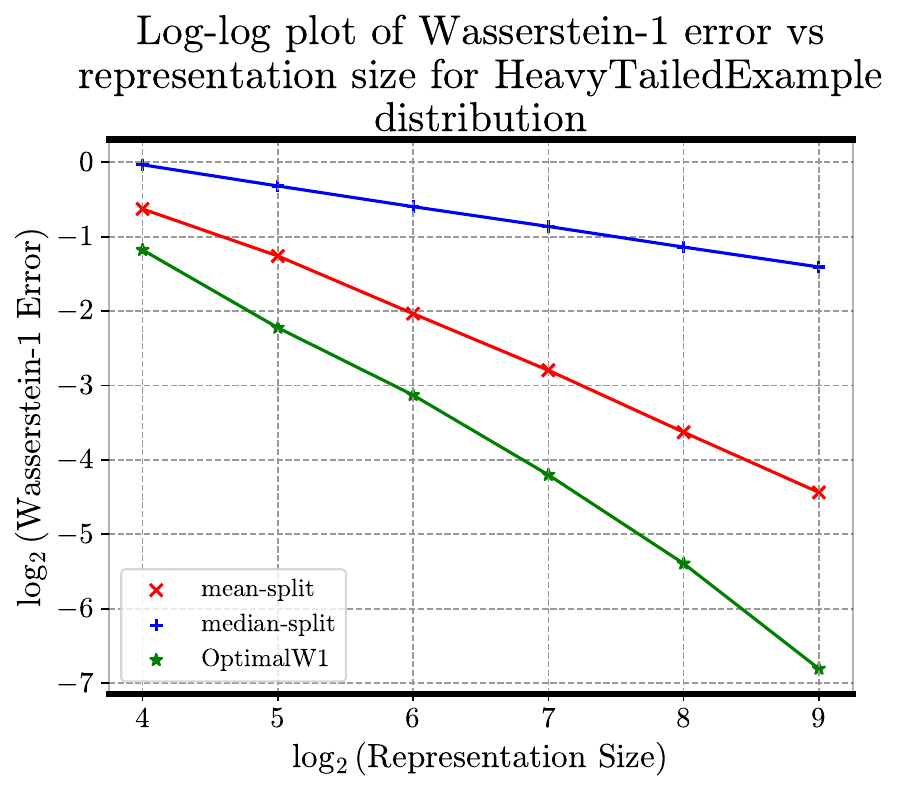}
    \end{subfigure} 
    \begin{subfigure}{0.49\textwidth}
        \centering
        \includegraphics[width=\linewidth]{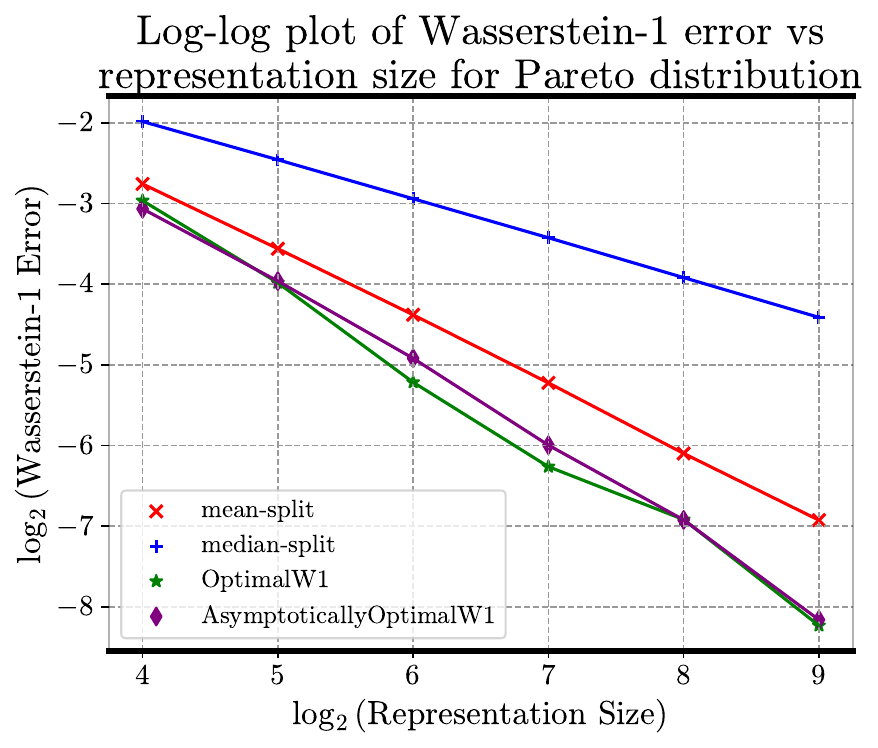}
    \end{subfigure} 
    \begin{subfigure}{0.49\textwidth}
        \centering
        \includegraphics[width=\linewidth]{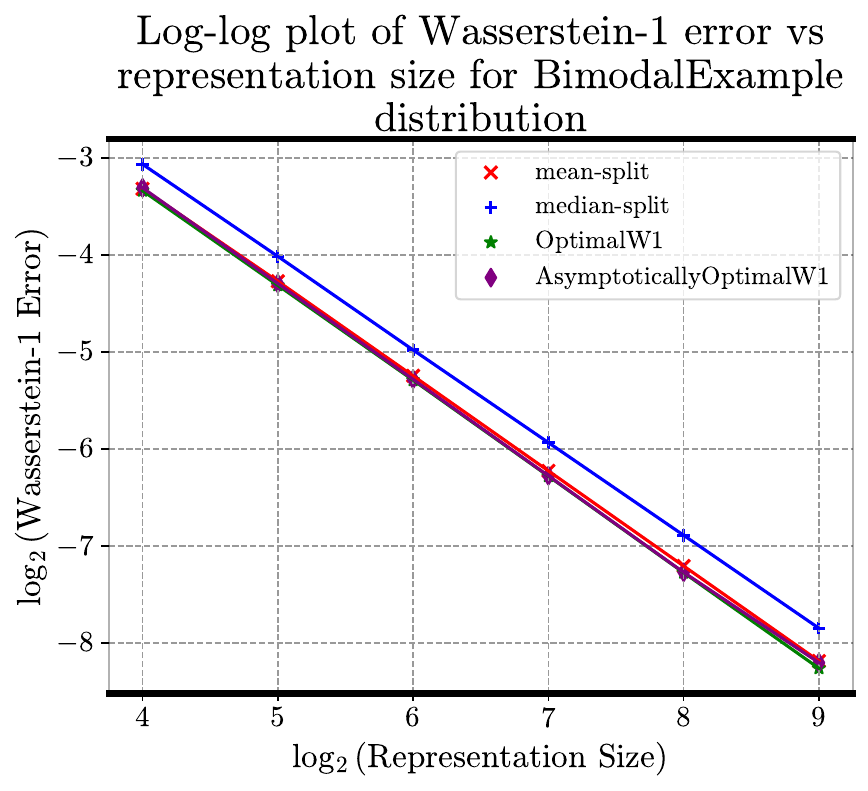}
    \end{subfigure} 

\caption{Log-log plots of Wasserstein-1 error vs representation size for mean-split, median-split, optimal and asymptotically optimal representations for the distributions in Table~\ref{tab:distributions}.}    \label{fig:log_wass_vs_repsize}
\end{figure}

Figure~\ref{fig:wass_vs_repsize} and Figure~\ref{fig:log_wass_vs_repsize} shows the plots of Wasserstein-1 error vs representation size for the Gaussian, Exponential, Pareto, heavy-tailed and bimodal example distributions as defined in Table~\ref{tab:distributions} as well as the equivalent log-log plots.
The corresponding log-log plots are shown below which demonstrate the optimal scaling for the mean-split algorithm.
Figures~\ref{fig:wass_vs_repsize} and~\ref{fig:log_wass_vs_repsize} show that, as expected, the optimal representations produce the lowest Wasserstein-1 errors, with the asymptotically optimal representation obtaining similar accuracies for the Gaussian, Exponential and Pareto distributions where the conditions outlined in Equation~\eqref{e.asymp-opt-assumptions} are fulfilled.
One can observe that the mean-split algorithm is often relatively close to optimal for all distributions.
The asymptotically optimal representation is omitted for the heavy-tailed example distribution as the conditions outlined in Equation~\eqref{e.asymp-opt-assumptions} are not satisfied.

\subsection{Performing arithmetic operations between distributions}
An issue that arises from performing arithmetic operations on a discrete representation of probability distributions is that the arithmetic operations between said distributions lead to an exponentially large number of atoms. For example, if we were to add two distributions with a representation size of $N$, then this would yield $N^2$ terms. In general, the sum of $n$ representations yields (worst-case) a representation of $N^n$ distinct atoms.
To deal with this curse of dimensionality we compress the distribution after every arithmetic operation using a discrete version of the quantization algorithms discussed in this article from $N^2$ to $N$. In this article, we use the same quantization algorithm for both constructing and compressing discrete representations of distributions. In other words, to the mean-split representation we apply the mean-split algorithm to perform the compression. This compression approach allows us to keep the representation size fixed. Some care is needed when using the quantile-based algorithms (such as the median-split and the asymptotically-optimal) as a compression algorithm. Indeed it is easy to show that given a discrete distribution $\mu = p \delta_x + (1-p) \delta_y$ we have $\mu(\Omega_-) \in \{0,1\}$. For this reason we linearly interpolate between atoms to avoid this issue.

As a concrete example, we explain how the compression of a $k$-fold convolution is carried out. Suppose $T(\mu,n)$ is defined using $f(\mu) = \bar{\mu}$. Let 
\begin{equation}
    \CC : \CM \mapsto \{ \mu \in \CM : \#(\supp(\mu)) = 2^n\},
\end{equation}
denote a compression algorithm to be defined below. Given $\mu_i, i=1,2,...,k$ and their discrete representations $\mu_i^{(n)}, i = 1,2,\dots,k$ define 
$\nu_j^{(n)}, j=1,2,\dots,k-1$ by
\[
\nu_{1}^{(n)} = T(\mu_1^{(n)}*\mu_2^{(n)},n),\ \nu_{i+1}^{(n)} = T(\mu_{i+2}^{(n)}*\nu_{i}^{(n)},n),\ i = 1,2,\dots,k-2.
\]
The compression of the $k$-fold convolution is then defined by \[
\CC(\mu_1^{(n)} *...*\mu_k^{(n)}) = \nu_{k-1}^{(n)}.
\]
We note that the compression procedure defined above is in general not associative and the result may depend on the order in which the operations are performed.

\begin{figure}[H]
    \centering
    \begin{subfigure}{0.49\textwidth}
        \centering
        \includegraphics[width=\linewidth]{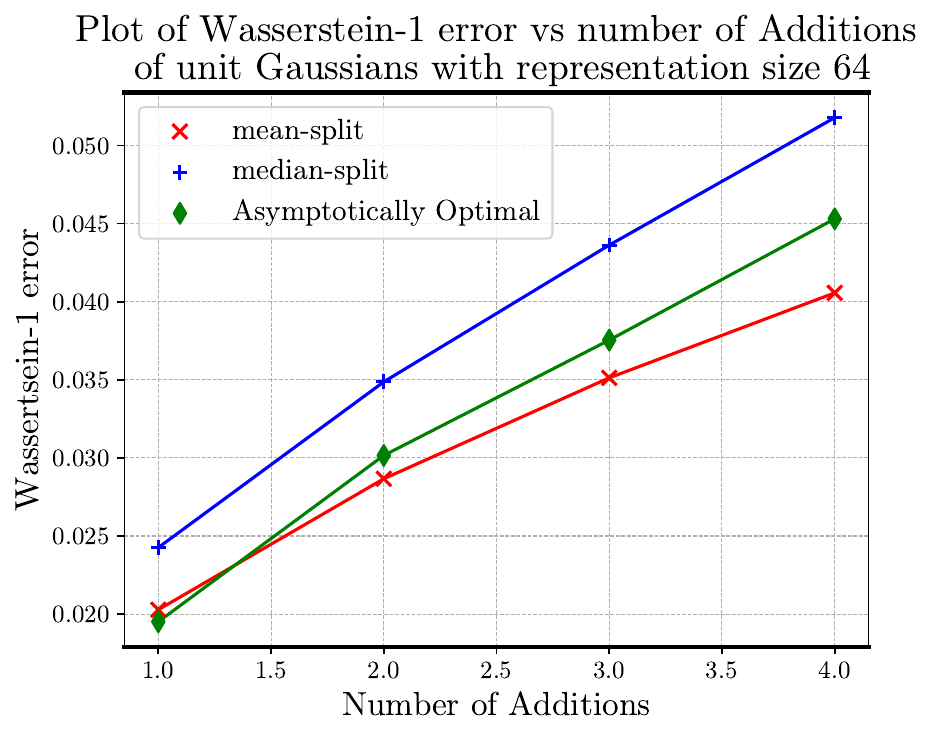}
    \end{subfigure} 
    \begin{subfigure}{0.49\textwidth}
        \centering
        \includegraphics[width=\linewidth]{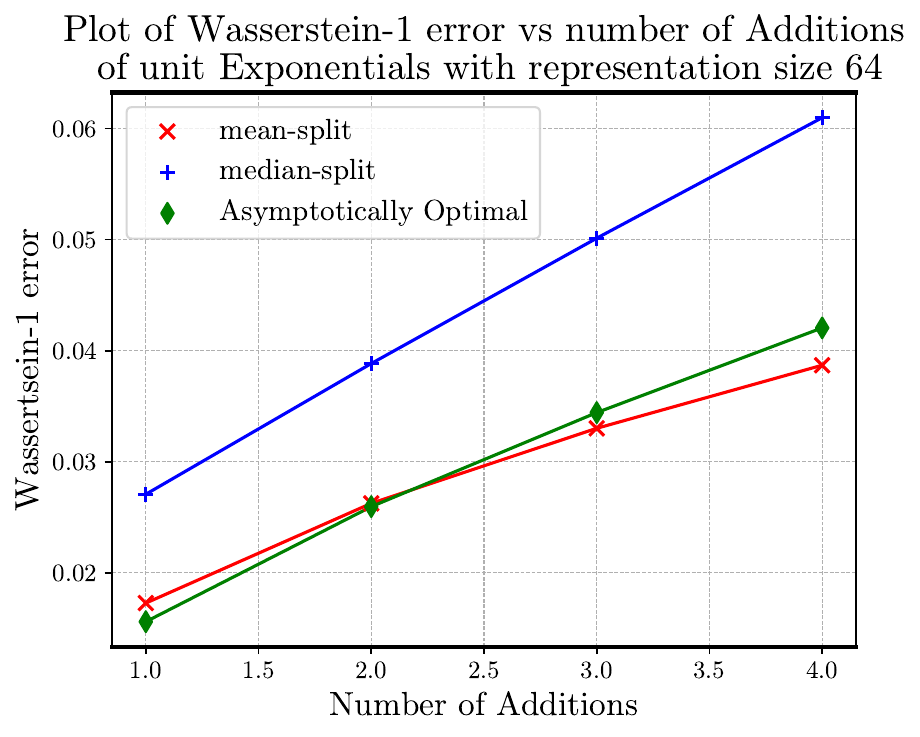}
    \end{subfigure}
    \begin{subfigure}{0.49\textwidth}
        \centering
        \includegraphics[width=\linewidth]{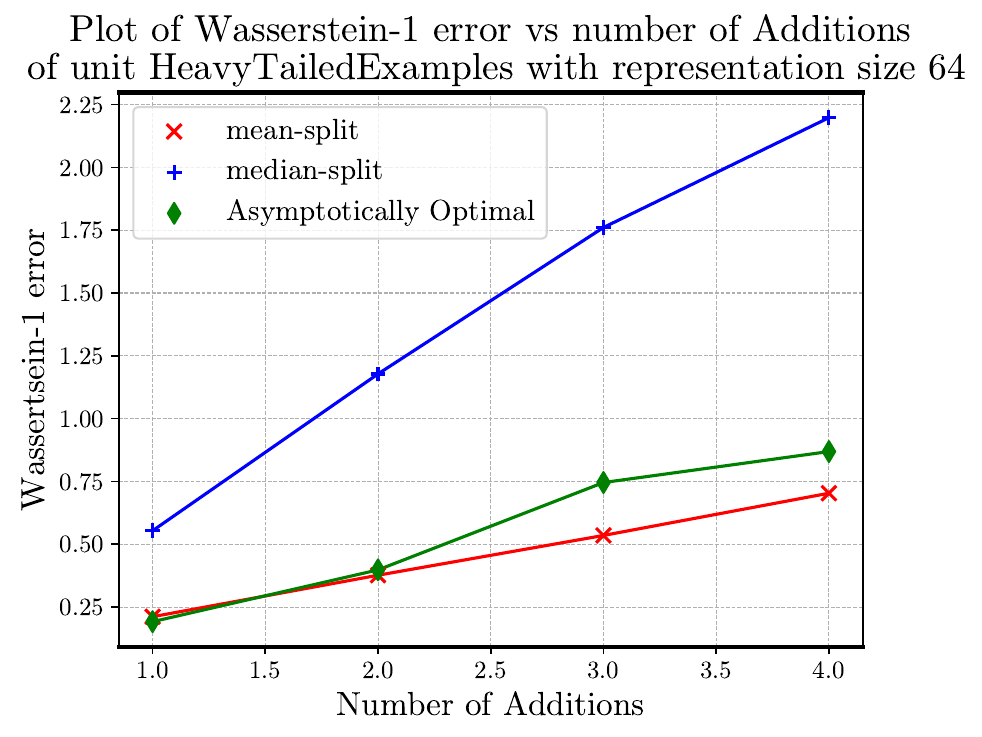}
    \end{subfigure} 
    \begin{subfigure}{0.49\textwidth}
        \centering
        \includegraphics[width=\linewidth]{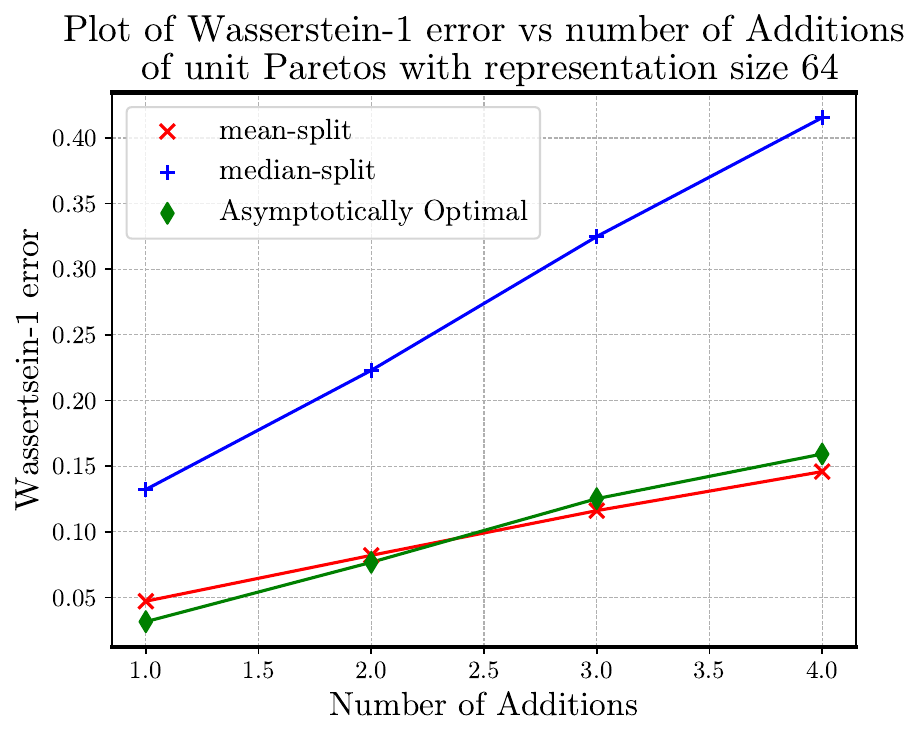}
    \end{subfigure}
    \begin{subfigure}{0.49\textwidth}
        \centering
        \includegraphics[width=\linewidth]{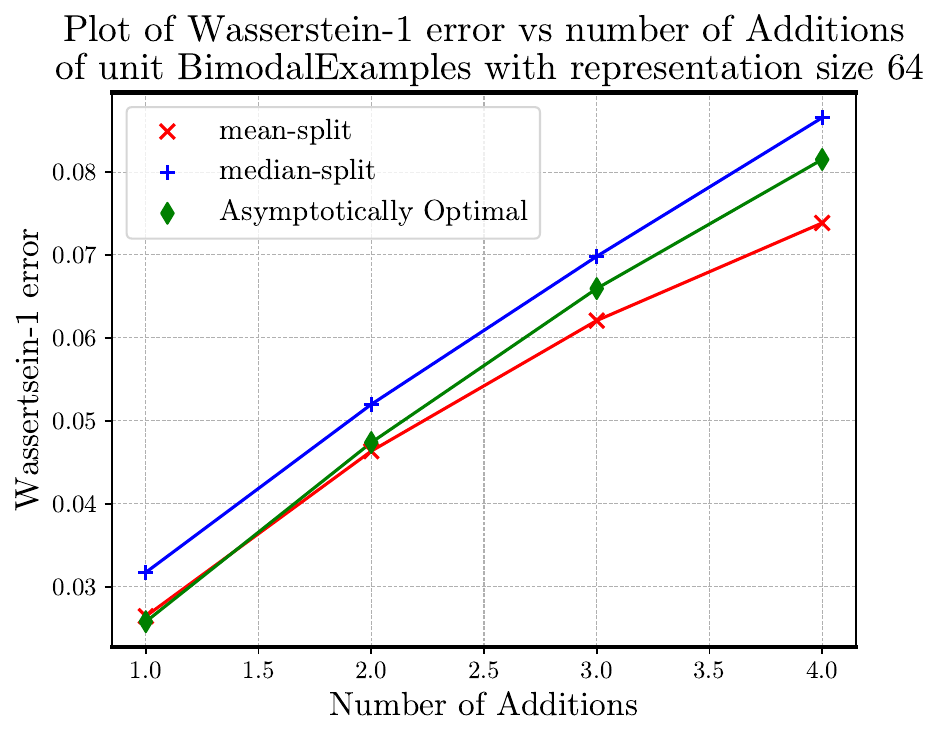}
    \end{subfigure}
    \caption{Wasserstein-1 error vs the number of additions for mean-split, median-split and 
    asymptotically-optimal representations for the distributions in Table~\ref{tab:distributions} 
    with representation size 64.}    \label{fig:wass_vs_num_additions}
\end{figure}

\begin{figure}[H]
    \centering
    \begin{subfigure}{0.49\textwidth}
        \centering
        \includegraphics[width=\linewidth]{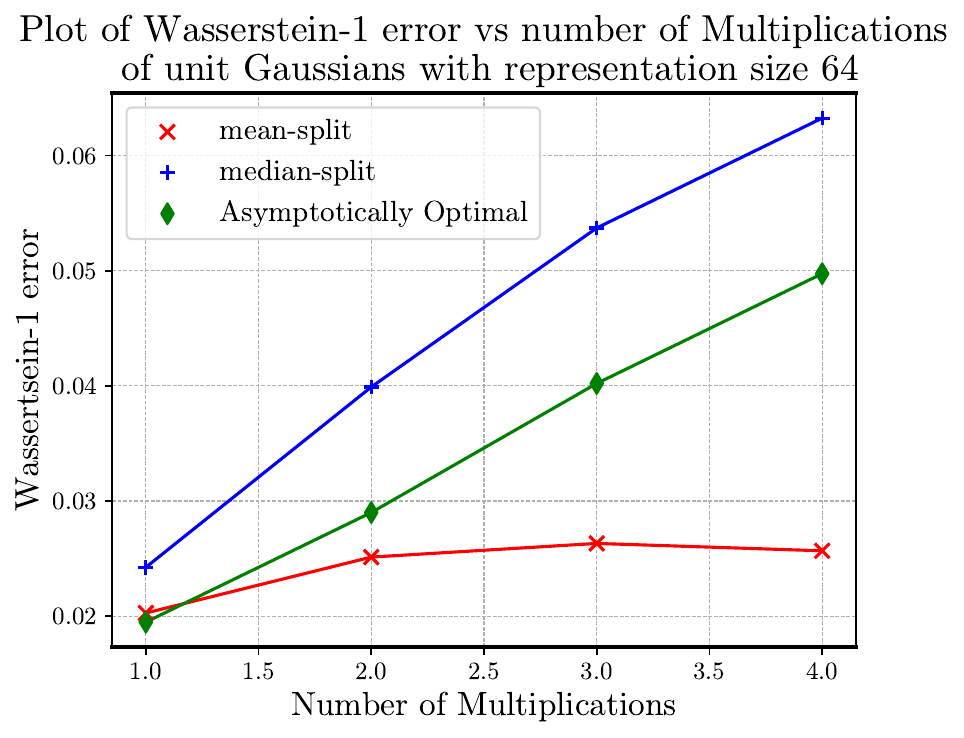}
    \end{subfigure} 
    \begin{subfigure}{0.49\textwidth}
        \centering
        \includegraphics[width=\linewidth]{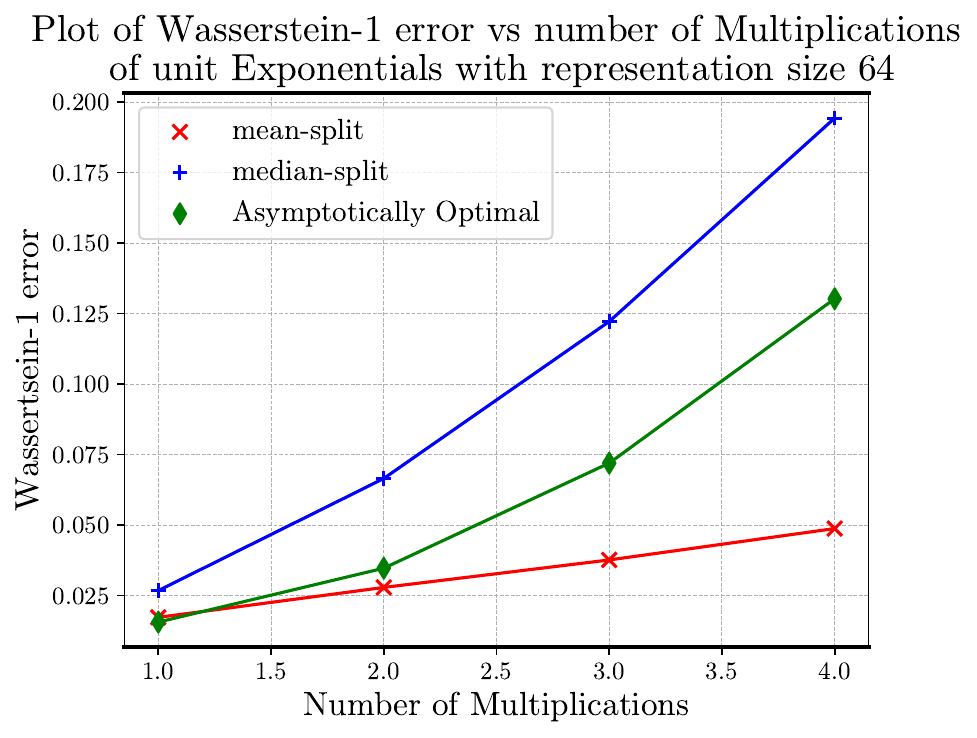}
    \end{subfigure}
        \begin{subfigure}{0.49\textwidth}
        \centering
        \includegraphics[width=\linewidth]{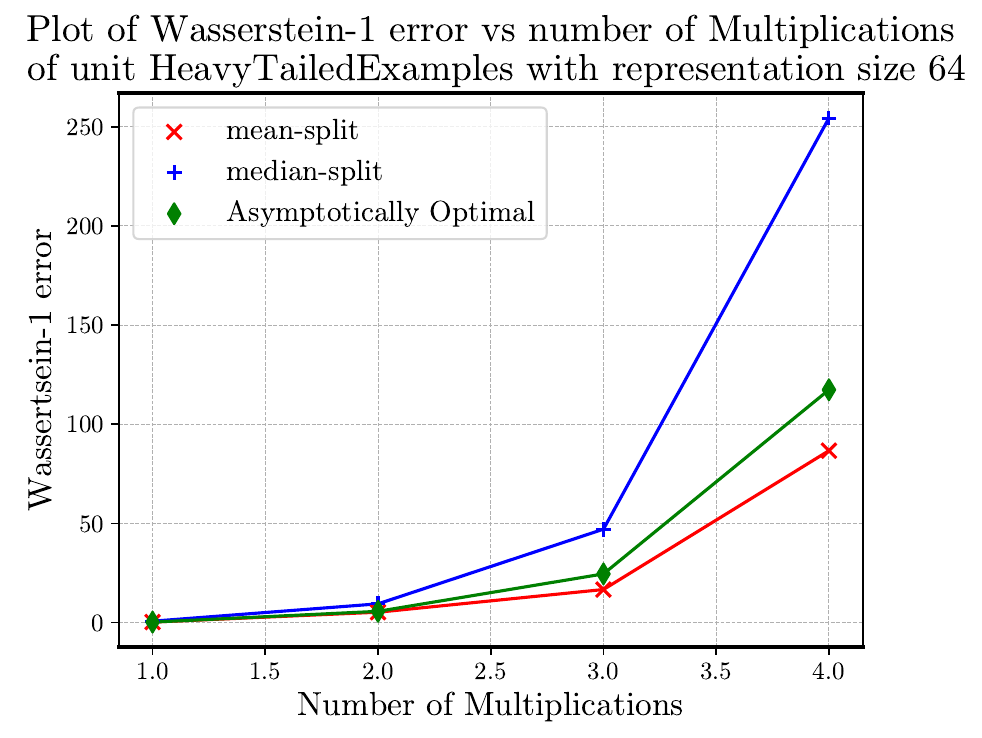}
    \end{subfigure} 
    \begin{subfigure}{0.49\textwidth}
        \centering
        \includegraphics[width=\linewidth]{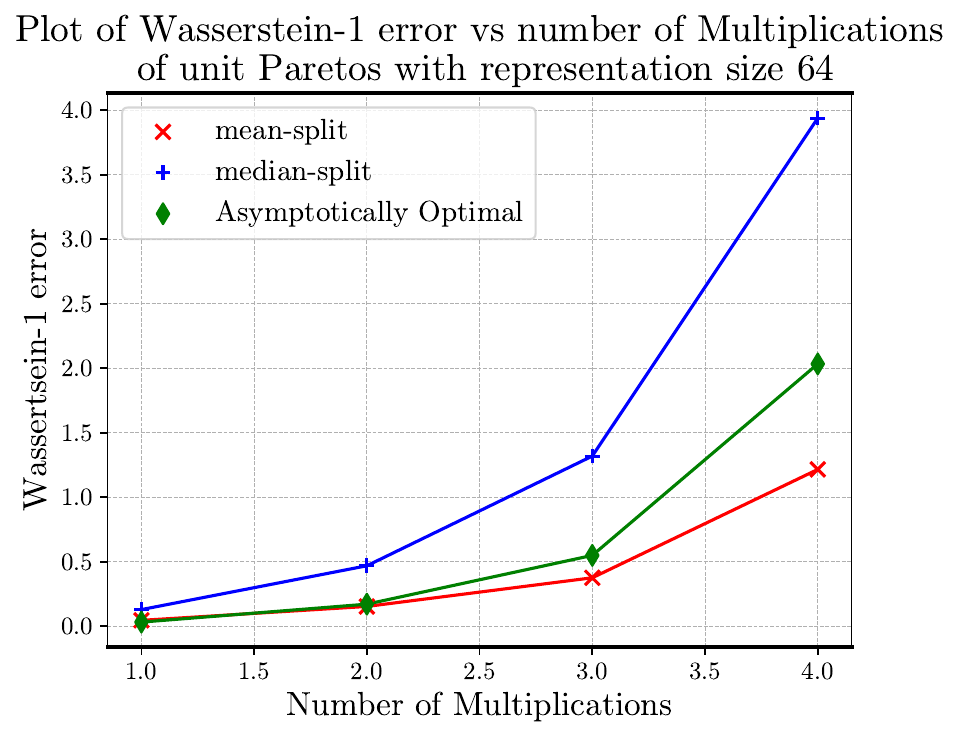}
    \end{subfigure}
    \begin{subfigure}{0.49\textwidth}
        \centering
        \includegraphics[width=\linewidth]{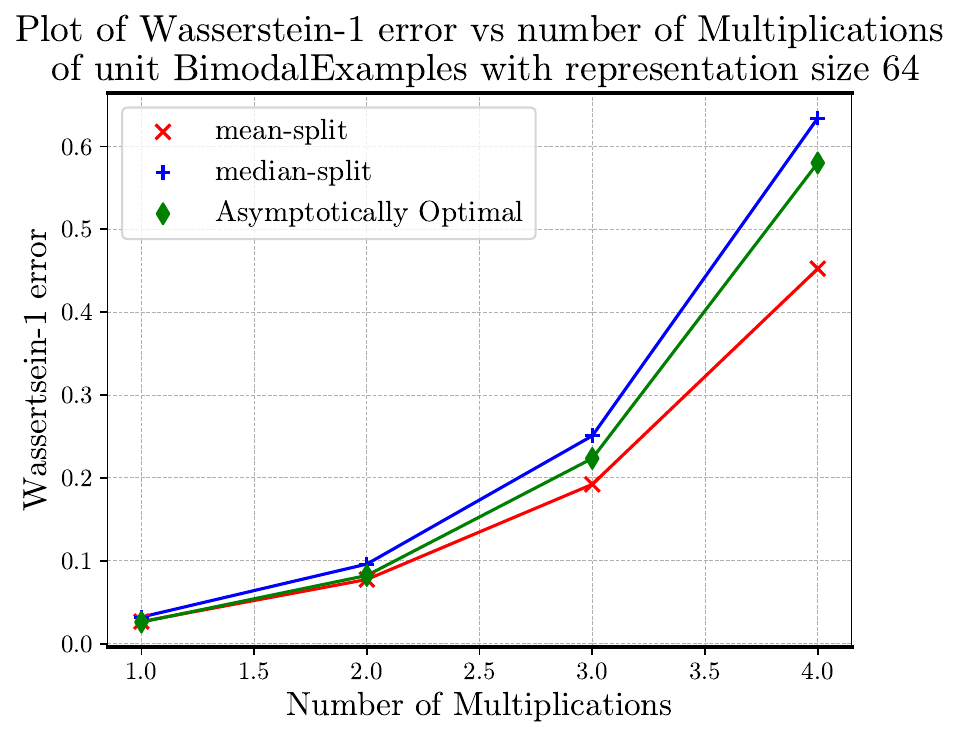}
    \end{subfigure}
    \caption{Wasserstein-1 error vs the number of multiplications for mean-split, median-split and 
    asymptotically-optimal representations for the distributions in Table~\ref{tab:distributions} 
    with representation size 64.}
\label{fig:wass_vs_num_multiplications}
\end{figure}

Before presenting the plots, we also comment on the error introduced by the compression step. The Wasserstein-1 error between the convolution \( \mu_1^{(n)} * \mu_2^{(n)} \) and its compressed version \( T(\mu_1^{(n)} * \mu_2^{(n)}, n) \) can be bounded from above using \eqref{e.discrete_upper_bound} as follows. Let $M_i = \max(\supp(\mu_i^{(n)})),\ m_i = \min(\supp(\mu_i^{(n)}))$ then
\begin{equation}\label{e.compression-bound}
    \begin{split}
        W_1(\mu_1^{(n)} * \mu_2^{(n)}, T(\mu_1^{(n)} * \mu_2^{(n)}, n)) 
        &\leq   \frac{M_1-m_1+M_2-m_2 }{2^{n+1}}.
    \end{split}
\end{equation}
Figures~\ref{fig:wass_vs_num_additions} and~\ref{fig:wass_vs_num_multiplications} show how the Wasserstein-1 error varies as a function of the number of additions/multiplications of all the distributions in Table~\ref{tab:distributions}, where the distributions are kept at a constant representation size of 64 using the compression approach with a range of discrete representation algorithms for the compression. We omit the optimal representation from the following plots as its dependence on general optimization methods and numerical solvers make it too computationally expensive for repeated distributional arithmetic and compression operations.

Figures~\ref{fig:wass_vs_num_additions} and~\ref{fig:wass_vs_num_multiplications} show that more accurately representing two distributions does not mean that you will better represent the sum/product of said distributions. For instance, Figure~\ref{fig:wass_vs_num_additions} shows that the sum of two unit Gaussian distributions computed with the mean-split compression algorithm has a lower Wasserstein-1 error than the same result computed using the asymptotically-optimal representation. This is despite the fact that a single unit Gaussian distribution is more accurately represented (with respect to Wasserstein-1) by the asymptotically optimal representation.
Figures~\ref{fig:wass_vs_num_additions} and~\ref{fig:wass_vs_num_multiplications} demonstrate that the mean-split algorithm performs better for addition and multiplication than the asymptotically-optimal and median-split algorithms for all distributions and arithmetic operations and is the best performing algorithm overall for repeated arithmetic and compression operations.

\subsection{Comparison with Monte Carlo}
As discussed in Section~\ref{s.intro}, Monte Carlo (MC) provides a popular and simple method to approximate distributions. Here we perform a simple comparison of MC and the mean-split algorithm. We want to determine the number of MC samples required such that MC and the mean-split algorithm produce the same Wasserstein-1 error. However, due to the stochastic nature of MC, we must instead calculate the mean Wasserstein-1 error or a large number of MC simulations. The distribution of the Wasserstein-1 error is known in the asymptotic limit \cite{MR1698999}. Consider a distribution $\mu$ with CDF $F_\mu$. Let $W_{1}^{(n)}$ denote the Wasserstein-1 error between $\mu$ and the empirical measure $\mu^{(n)}$ generated by $n$ i.i.d samples from $\mu$. Then due to \cite[Theorem 2.1]{MR1698999} the expected value of the Wasserstein-1 error is asymptotically equal to
\begin{equation*}
    \sqrt{n}\ \mathbb{E}\left(W_{1}^{\left(n\right)}\right)=(1+o(1))\sqrt{\frac{2}{\pi}}\int_{\BR}\sqrt{F\left(x\right)\left(1-F\left(x\right)\right)}\id x.
\end{equation*}
For the unit Exponential distribution, we see that $\sqrt{n}\ \mathbb{E}\left(W_{1}^{(n)}\right)=(1+o(1))\sqrt{\pi/2}$. We can calculate that a mean-split representation with representation size 256 yields a Wasserstein-1 error of 0.0044371522. Hence, we may deduce that we would (by ignoring the $1+o(1)$ term) require 82,197 MC samples to produce an error lower than this value in expectation. Similarly the equivalent MC count of the unit Gaussian for 256 is 61,341. In fact, due to the $1/\sqrt{n}$ vs $1/n$ convergence of Wasserstein-1 errors for MC and the mean-split representation, we require quadratically more MC samples to achieve the same accuracy as the mean-split representation.

This comparison is especially important when considering distributions whose distribution or density functions are not known and must be treated computationally. In these cases, MC provides the most ubiquitous method for computing such distributions. However, the compression method allows for an alternative approach. We note that while the quantization algorithm presented in this article does require knowledge of the distribution (e.g., the ability to evaluate the CDF or compute conditional means), the compression step does not. Indeed, as noted in Remark~\ref{r.mean-split-remarks}, the algorithm applies equally well to discrete distributions. Thus, once discrete representations have been obtained, whether through the quantization algorithm, from empirical data, or from any other source, arithmetic operations and subsequent compression can be carried out without any knowledge of the underlying continuous distribution. Consider the product of four mutually-independent unit Gaussian random variables. One can determine that application of the mean-split algorithm with a representation size 512 is equivalent to 165,910 MC samples. This demonstrates how the mean-split algorithm provides an efficient and deterministic alternative to approximating such distributions.

Although we decided to compare the mean value of the distribution of Wasserstein-1 errors of the Monte Carlo simulations to that of the mean-split method, we note that this may not be always be the most appropriate method to use.  Another strategy would involve utilizing the asymptotic distribution of $W_1^{(n)}$, for example, to construct confidence intervals that ensure convergence at a specified confidence level. For example, if one was interested in the average accuracy across many MC simulations, then using the mean comparison method would be appropriate. If however, one wanted to find the minimum MC count that guarantees a certain degree of convergence, then the equivalent representation size for the mean-split method would be best found using the quantile approach.

\section{Discussion and summary}\label{s.discussion}
This article introduced a simple divide-and-conquer algorithm for approximating one-dimensional distributions. We provided a general and easily computable upper bound on the Wasserstein distance between the target distribution and its approximation, under minimal assumptions. The effectiveness of this bound was demonstrated in Sections~\ref{ss.examples} and~\ref{ss.corollaries}. We also analyzed the stability of the representations with respect to arithmetic operations. Our numerical results suggest that the mean-split algorithm outperforms the quantile-based algorithm in certain scenarios. 

What follows are observations and remarks that, while relevant, do not naturally fit within the main body of the article.
\subsection{Comments on the computational complexity and robustness of the algorithm}
To provide a more complete picture of the effectiveness of the introduced algorithm from a practical standpoint, we address the computational complexity of the algorithms discussed in this article. 

We first address the case when $\mu$ is continuous as the discrete case depends on the number of atoms. We shall assume that the given split function $f$ is computable in constant time. This is likely not true in general and depends on the distribution. In the context of the split functions discussed in this article (i.e the mean and median) it depends on how fast one can evaluate either of the functions 
\begin{equation*}
\begin{split}
&\int_{-\infty}^x t \id \mu(t) ,\ x \in \BR,\\
&F^{-1}(p),\ p\in [0,1].
\end{split}    
\end{equation*}
In practice, this relies on efficient implementations of special functions such as the incomplete Gamma or Beta functions, see Press et. al. \cite{press2007numerical} for an extensive treatment on this topic. Nevertheless, under this assumption the computational complexity of $T(n,\mu)$, denoted $\SFC(n)$, satisfies the recursion 
\[
\SFC(n) = 2 \SFC(n-1) +\CO(1).
\]
Hence, the complexity satisfies $\SFC(n) = \CO(2^n)$. In other words, the algorithm is linear in the representation size $N=2^n$. Under a similar assumption, the asymptotically-optimal representation also has computational complexity that is linear in the representation size. 

For a discrete distribution, 
\[
\mu = \sum_{i=1}^m p_i \delta_{x_i},\ -\infty < x_1<x_2<...<x_m<\infty,
\] 
we assume that the evaluation of $f(\mu)$ is linear in $m$. This is certainly true for the mean-split function. The computational complexity of the mean-split method is $\CO(mn)$. Indeed, for any fixed $\mu$, the computational complexity, $\SFC(n,m)$ satisfies the following recursion. For any $n\geq 0$ and $m\geq 2^n$, there exist $m_1,\ m_2\geq 1,\  m_1+m_2 = m$ (depending on $\mu$) such that 
\[
\SFC(n,m) = \SFC(n-1,m_1)+\SFC(n-1,m_2) +\CO(m), \SFC(0,m) = \CO(m).
\]
From this it follows by induction that for any $k = 1,2,\dots,n$:
\[
\SFC(n,m) = \sum_{i=1}^{2^k} \SFC(n-k,m_i) +\CO(km),
\]
where $m_i\geq 1, \forall i$ and $\sum_{i=1}^{2^k} m_i = m$. Hence 
\[
\SFC(n,m) = \CO(nm)+\sum_{i=1}^{2^n} \SFC(0,m_i)  =  \CO(nm)+\sum_{i=1}^{2^n}  \CO(m_i ) = \CO(nm),
\]
as claimed. On the other hand, the computational complexity of the asymptotically optimal method is \( \CO(m \cdot 2^n) \). In the context of compressing a distribution with \( N^2 \) atoms back to \( N \) atoms, the mean-split method has computational complexity \( \CO(N^2 \log N) \), whereas the asymptotically optimal method has complexity \( \CO(N^3) \).

Finally, the mean-split algorithm demonstrates a form of robustness when applied to discrete distributions. As noted in Remark~\ref{r.lower_bound}, it preserves the mean of the original distribution. Moreover, by Definition~\ref{def.algorithm}, we have:
\[
T(\delta_x,n) = \delta_x,\ x \in \BR,\ n \geq 0,
\]
then the mean-split algorithm satisfies 
\begin{equation}\label{e.mean-split-identity-property}
     T(\mu,n) = \mu, \forall n\geq \#(\supp(\mu))-1.
\end{equation}
While we shall not give a formal proof (as we do not need it in the paper) we hope that the following reasoning provides a compelling enough argument. 
The set of discrete probability measures with fixed number of atoms can be mapped into the set of all finite binary trees (i.e., a tree in which each node has at most two children) due to the recursive structure of the algorithm. Thus, the key question becomes determining the maximum height that these trees can take. Note that since $\#(\supp(\mu_+)) + \#(\supp(\mu_-)) = \#(\supp(\mu)) = m$, the worst-case scenario occurs when either $\#(\supp(\mu_+)) = 1$ or $\#(\supp(\mu_-)) = 1$. In such a case, the algorithm effectively returns to the situation described just above Equation~\eqref{e.mean-split-identity-property}, but with \( n \) decreased by 1. Iterating this process, we see that at iteration \( n - 1 \), the worst-case scenario leaves us with a distribution \( \nu \) with \( \#(\supp(\nu)) = 2 \), from which it is easy to verify that \( T(1, \nu) = \nu \). For any other distribution the height of the tree will be smaller.

This is not true for median-based split functions. For example, when applied to a distribution with only $2$ atoms, $\mu = p \delta_x + (1-p) \delta_y$ it will either put all mass on $x$ or $y$ depending on $p \in (0,1).$ Interestingly, the median-split algorithm, which corresponds to the greedy choice in the \( W_1 \) setting, appears to perform considerably worse than the mean-split algorithm in all considered examples. Analogously, the mean-split algorithm would correspond to the greedy choice were we to conduct the analysis in terms of the \( W_2 \) metric.
\subsection{Comments on the numerical experiments}
The numerical experiments indicate that the mean-split algorithm outperforms the other algorithms 
in the context of performing arithmetic operations between distributions. However, additional 
results and a more comprehensive numerical study are needed before drawing definitive conclusions. 
To advance this line of research, it would be valuable to extend the study to real-world 
applications, such as numerical approximations of SDEs. For example, the 
Euler-Maruyama scheme for approximating 
solutions of SDEs involves repeated arithmetic and functional operations 
on random variables at each time step. The arithmetic and compression 
methodology outlined in Section~\ref{s.numerical} could be used to 
propagate quantized distributions through these operations directly, 
providing a deterministic alternative to Monte Carlo-based methods.

Furthermore, it is important to assess the performance of compression algorithms within this context. 

While it would be very interesting to obtain theoretical results relating to the behavior of different representations with regards to arithmetic operations, it seems hard to deduce anything satisfactory. A prominent reason is that it is hard to find natural assumptions. Indeed, consider the following informal argument that highlights the issue. Let $\mu,\ \nu \in \CM_1$ and suppose (for this arguments sake) that $\Phi : \BR \times \BR \to \BR$ is Lipschitz continuous on $\supp(\mu)\times \supp(\nu)$ (this could likely be replaced by a combination of integrability and smoothness assumptions of the partial derivatives of $\Phi$). Then the Wasserstein-1 distance between $ (\mu \otimes \nu) \circ \Phi^{-1}$ and $(\mu^{(n)} \otimes \nu^{(n)}) \circ \Phi^{-1}$ is bounded by the distance between $\mu \otimes \nu$ and $\mu^{(n)} \otimes \nu^{(n)}$. Indeed, if $(X,Y) \sim \mu\otimes\nu$ and $(X_n,Y_n) \sim \mu^{(n)} \otimes\nu^{(n)}$, then for any coupling between the two pairs we have by the Lipschitz continuity of $\Phi$
\begin{equation*}
    \BE [|\Phi(X,Y)-\Phi(X_n,Y_n)|]
    \leq 
    C(\Phi) 
    (\BE[|X-X_n| +|Y-Y_n|]),
\end{equation*}
or in other words
\begin{equation}\label{e.generic-lipschitz-bound}
    W_1((\mu \otimes \nu) \circ \Phi^{-1},(\mu^{(n)} \otimes \nu^{(n)}) \circ \Phi^{-1})
    \leq 
    C(\Phi) 
    ( W_1(\mu,\mu^{(n)}) +W_1(\nu,\nu^{(n)})).
\end{equation}
Clearly this type of argument holds for any discrete representation and without the appropriate assumptions it is impossible to see any differences between the representations from a theoretical point of view. Thus, to better understand the effect of any such operation we have chosen to resort to a numerical study.

We remark that the only rigorous statement we can make in this context is that if the discrete 
approximation does not capture the mean exactly, the Wasserstein distance will grow as 
$\Theta(k)$ with $k$ additions. To be more precise, suppose that $\mu \in \CM_1$ and $\mu^{(n)}$ 
is obtained via the some split function $f$ with $\bar{\mu}\neq \bar{\mu}^{(n)}$. Then by \eqref{e.generic-lipschitz-bound} we can deduce the upper bound (using $\Phi(x,y) = x+y, C(\Phi)=1$):
\begin{equation*}
W_1
\left(
\underbrace{\mu *  \dots*  \mu}_{k \text{ times}}
,\underbrace{\mu^{(n)} *  \dots*  \mu^{(n)}}_{k \text{ times}}
\right)
\leq k W_1(\mu,\mu^{(n)}),
\end{equation*}
while the matching lower bound follows from the Lipschitz continuity of the mean with 
respect to $W_1$:
\[
k 
\left|
    \int t \id \mu (t) - \int t \id \mu^{(n)}(t)
\right|
\leq 
W_1
\left(
\underbrace{\mu *  \dots*  \mu}_{k \text{ times}}
,\underbrace{\mu^{(n)} *  \dots*  \mu^{(n)}}_{k \text{ times}}
\right).
\]
This phenomenon occurs when \( \mu \) is the Exponential distribution and \( \mu^{(n)} \) corresponds to the optimal quantization, as described in \cite[p.~70, Example~5.7]{MR1764176}.
\subsection{Comparing deterministic arithmetic on distributions and Monte Carlo simulations}
The final topic is a comparison between the deterministic approach and the Monte Carlo method. For the sake of this discussion, consider the following simple example. Let $F : \CM \times \CM \mapsto \CM$ be the map defined by
\[
F(\mu_1, \mu_2) = \mu_1 * \mu_2,
\]
which is Lipschitz continuous with respect to the \( W_1 \) metric, satisfying \( \|F\|_{\text{Lip}} \leq 1 \). Suppose that \( \mu_1 \) and \( \mu_2 \) satisfy the assumptions of Theorem~\ref{thm.poly-tails} with \( \alpha > 2 \). This choice ensures that both the Monte Carlo and deterministic approaches achieve their theoretically optimal rates of convergence.

In the standard Monte Carlo approach, one constructs the empirical measure
\[
F^{(N)} = \frac{1}{N} \sum_{i=1}^N \delta_{X_1^i + X_2^i},
\]
where \( (X_1^i, X_2^i) \sim \mu_1 \times \mu_2 \), for \( i = 1, 2, \dots, N \), are i.i.d.\ samples. By \cite[Theorem 1.1]{MR1698999} we have the convergence in distribution
\begin{equation}\label{e.MC-convergence}    
\sqrt{N} \, W_1(F, F^{(N)}) \xrightarrow{d} \int_0^1 |B(t)| \, \mathrm{d}Q(t),
\end{equation}
where \( Q : [0,1] \to \mathbb{R} \) is the quantile function of the distribution \( \mu_1 * \mu_2 \), and \( B(t) \) denotes a standard Brownian bridge. Heuristically, this result suggests that \( W_1(F, F^{(N)}) \approx \frac{C}{\sqrt{N}} \), where \( C \) is a random constant. However, to quantify how large this constant might be, one would need to know the distribution of \( F \) and its quantile function, which is typically unavailable in practice.

In comparison, the deterministic approach allows for the construction of a non-compressed approximation via
\[
F^{(n)} = \mu_1^{(n)} * \mu_2^{(n)} = \sum_{(x_1, x_2)} \mu_1^{(n)}(\{x_1\}) \mu_2^{(n)}(\{x_2\}) \delta_{x_1 + x_2},
\]
or a compressed approximation given by
\[
F^{(n),c} = T(\mu_1^{(n)} * \mu_2^{(n)}, n).
\]
Due to the curse of dimensionality, the non-compressed approximation, while accurate, quickly becomes computationally infeasible as the number of terms in the sum grows exponentially with respect to the number of convolutions. In practice, it is therefore more efficient to use the compressed approximation \( F^{(n),c} \).
From~\eqref{e.compression-bound}, we obtain the following upper bound on the error between \( F \) and \( F^{(n),c} \):
\begin{align*}
    W_1(F^{(n),c},F) &\leq W_1(T(\mu_1^{(n)} * \mu_2^{(n)},n), \mu_1^{(n)} * \mu_2^{(n)})+ W_1(\mu_1^{(n)} * \mu_2^{(n)}, \mu_1 * \mu_2) \\
    & \leq W_1(T(\mu_1^{(n)} * \mu_2^{(n)},n), \mu_1^{(n)} * \mu_2^{(n)})+ W_1(\mu_1^{(n)} , \mu_1 )+W_1( \mu_2^{(n)},  \mu_2) \\
    &= \frac{1}{2}\frac{m_n-l_n}{2^{n}}+W_1(\mu_1^{(n)} , \mu_1 )+W_1 ( \mu_2^{(n)},  \mu_2),
\end{align*}
where 
\begin{align*}    
m_n &= 
\max
    \left(
        \supp\left(\mu_1^{(n)}\right)
    \right)
+
\max
    \left(
        \supp\left(\mu_2^{(n)}\right)
    \right),\\
l_n &= 
\min
    \left(
        \supp\left(\mu_1^{(n)}\right)
    \right)
+
\min
    \left(
        \supp\left(\mu_2^{(n)}\right)
    \right).
\end{align*}
By the assumptions on $\mu_1,\ \mu_2$ we see that $W_1(\mu_1^{(n)} , \mu_1 )+W_1( \mu_2^{(n)},  \mu_2) = \Theta(2^{-n(1+o(1))})$ while the first follows a similar bound. Thus 
\[
 W_1(F^{(n),c},F) \leq C 2^{-n(1+o(1))},
\]
with a deterministic constant $C$.
Letting $N = 2^n$ denote the representation size and assuming that the compression step has computational complexity \( \CO(N^2 \log(N) ) \), the overall complexity of the deterministic approach is \( \CO(N^2 \log(N) )\) with a deterministic constant. In contrast, achieving a comparable error using the Monte Carlo method requires only \( \CO(N^2) \) computational effort, albeit with a random constant. In practice, ensuring convergence in the Monte Carlo approach typically requires constructing confidence intervals based on the asymptotic distribution in~\eqref{e.MC-convergence}. This, in turn, involves computing the empirical quantile function, which requires sorting the samples generated by the Monte Carlo method, thereby increasing the computational complexity to \( \CO(N^2 \log(N)) \). This also introduces additional computational overhead, as one must estimate the CDF of the random variable appearing in~\eqref{e.MC-convergence}. From this, an interesting recursive situation also arises: to validate the results of the Monte Carlo simulation one needs to carry out another Monte Carlo simulation.
\subsection{Open questions}
We end this discussion with a list of open questions and problems, some of which have already been mentioned. 
The numerical results of this article seems to indicate that, at least in certain scenarios, the mean-split algorithm performs better than the asymptotically-optimal one. This leads to the following questions. 
\begin{question}
    Assume $f(\mu) = \bar{\mu}$. For which distributions (other than the uniform distribution) is $T(\mu,n)$ asymptotically optimal in the sense of Graf and Luschgy \cite[Section 7 p.93]{MR1764176}?
\end{question}
Moreover, from Example~\ref{ex.pareto} we see that the mean-split function also performs better throughout the entire valid parameter space (when compared to the median-split function) which leads us to the following question:
\begin{question}
    Is there a split function that performs uniformly better over all finite mean distributions? In other words, does there exist an optimal split function?
\end{question}
Finally two questions which were posed in Remark~\ref{r.polytail-conclusions}.
\begin{question}
\noindent
    \begin{itemize}
        \item 
        Assume that $\mu \in \CM_2,\ f(\mu)=\bar{\mu}$. Is it the case that 
        \[
        \limsup_{n\to \infty} \frac{\log W_1(\mu^{(n)},\mu)}{n} = -\log 2?
        \]
        \item
        Assume that there exists a $t_0>0$ such that $\int_\BR \e^{tx} \id \mu(x) < \infty$ for all $t \in (0,t_0)$. Is it the case that 
        \[
        W_1(\mu^{(n)},\mu) = \CO(2^{-n})?
        \]
    \end{itemize}
\end{question}
Now we turn to mention some more open-ended problems.
A problem not addressed in this paper is how one approximates distributions in higher dimensions. Indeed the algorithm relies heavily on the ability to split the support of the distribution using only a single datum specified by the split function. This leads to the following open problem.
\begin{problem}
Generalize the algorithm to higher dimensions. Find good candidate split functions in higher dimensions. 
\end{problem}
We also note that while we are able to obtain bounds for finite mean distributions we are at this moment unable to produce anything meaningful in the absence of finite mean. This leads us to pose:
\begin{problem}
    Extend the analysis for fat-tailed distributions, meaning distributions without finite mean. 
\end{problem}
On the discrete side, we note that the upper bound \eqref{e.discrete_upper_bound} is far from sharp and improving it would greatly help in better understanding the compression step when performing distributional arithmetic. This leads to the following problem:
\begin{problem}
    Improve the bound \eqref{e.discrete_upper_bound}.
\end{problem}
Finally, a much more open-ended problem. It is clear that the entire procedure assumes that one is able to evaluate the CDF analytically to obtain good bounds. Hence we state the following:
\begin{problem}
    Find an efficient algorithm for approximating distributions whose PDFs or CDFs are not analytically expressible (such as the $\alpha$-stable distributions).
\end{problem}

\newpage
\printbibliography

\end{document}